\documentclass[noams,lp]{compositio}

\usepackage{amsbsy,amssymb,amscd,amsfonts,latexsym,amstext,delarray,
 amsmath,stackrel,lineno,tabularx,url,MnSymbol,tikz,cite}

\input xypic
\usepackage[pdfborder={0 0 0}]{hyperref}

\hypersetup{
    colorlinks = true,
    linkcolor = blue,
    anchorcolor = red,
    citecolor = green,
    filecolor = red,
    pagecolor = red,
    urlcolor = magenta
}

\newcommand*\patchAmsMathEnvironmentForLineno[1]{%
  \expandafter\let\csname old#1\expandafter\endcsname\csname #1\endcsname
  \expandafter\let\csname oldend#1\expandafter\endcsname\csname end#1\endcsname
  \renewenvironment{#1}%
     {\linenomath\csname old#1\endcsname}%
     {\csname oldend#1\endcsname\endlinenomath}}%
\newcommand*\patchBothAmsMathEnvironmentsForLineno[1]{%
  \patchAmsMathEnvironmentForLineno{#1}%
  \patchAmsMathEnvironmentForLineno{#1*}}%
\AtBeginDocument{%
\patchBothAmsMathEnvironmentsForLineno{equation}%
\patchBothAmsMathEnvironmentsForLineno{align}%
\patchBothAmsMathEnvironmentsForLineno{flalign}%
\patchBothAmsMathEnvironmentsForLineno{alignat}%
\patchBothAmsMathEnvironmentsForLineno{gather}%
\patchBothAmsMathEnvironmentsForLineno{multline}%
}
\newcommand*{\DashedArrow}[1][]{\mathbin{\tikz [baseline=-0.25ex,-latex, dashed,#1] \draw [#1] (0pt,0.5ex) -- (1.3em,0.5ex);}}%

\newtheorem{thm}{Theorem}[section] 
\newtheorem{defn}[thm]{Definition} 
\newtheorem{prop}[thm]{Proposition}

\newtheorem{lem}[thm]{Lemma}

\newtheorem{rem}[thm]{Remark}

\def\Aut{{\rm Aut}}

\def\End{{\rm End}}

\def\Hom{{\rm Hom}}

\def\sign{{\rm sign}}

\def\sign{\mathbf{S}}

\def\B{{\mathbb B}}
\def\C{{\mathbb C}}
\def\F{{\mathbb F}}

\def\N{{\mathbb N}}
\def\P{{\mathbb P}}
\def\Q{{\mathbb Q}}
\def\R{{\mathbb R}}
\def\Z{{\mathbb Z}}

\def\cF{{\mathcal F}}
\def\cJ{{\mathcal J}}

\def\cP{{\mathcal P}}

\def\qqq{\,,\quad~\forall}

\newcommand{\ie}{{\it i.e.\/}\ }
\newcommand{\eg}{{\it e.g.\/}\ }
\newcommand{\cf}{{\it cf.}}
\newcommand{\opcit}{{\it op.cit.\/}\ }

\def\mod{semimodule}

\def\Mods{Semimodules}
\def\mods{semimodules}
\def\modsub{subsemimodule}
\def\modsubs{subsemimodules}

\def\tensb{{\otimes_\B\sign}}
\def\tensid{{\otimes_\B{\rm id}_\sign}}

\def\dim{{\mbox{dim}}}

\def\Hom {{\mbox{Hom}}}

\def\End{{\mbox{End}}}

\def\ff{{\rm  Fr}}

\def\rmax{\R_+^{\rm max}}

\def\cF{\mathfrak F}

\def\sign{\mathbf{S}}

\def\Se{\frak{ Sets}}
\def\fin{\frak{ Fin}}

\def\Arc{\frak{Arc}}

\def\proj{\frak{ Proj}}

\def\bff{\F^{\rm pf}}

\def\rk{{\rm rk}}
\def\ruk{\underline{\rm rk}}
\begin{document}

\title{Projective geometry in characteristic one and the epicyclic category}
\author{Alain Connes}
\email{alain@connes.org}
\address{Coll\`ege de France,
3 rue d'Ulm, Paris F-75005 France\newline
I.H.E.S. and Ohio State University.}
\author{Caterina Consani}
\email{kc@math.jhu.edu}
\address{Department of Mathematics, The Johns Hopkins
University\newline Baltimore, MD 21218 USA.}
%
\classification{19D55, 12K10, 51E26, 20N20.}
\keywords{Characteristic one, Cyclic category, Semifields, Projective geometry, Hyperstructures.}
\thanks{The second author is partially supported by the NSF grant DMS 1069218
and would like to thank the Coll\`ege de France for some financial support.
}

\begin{abstract}
We show that the cyclic and epicyclic categories which play a key role in the  encoding of cyclic 
homology and the lambda operations, are obtained from projective geometry in characteristic one
over the infinite semifield of ``max-plus integers'' $\Z_{\rm max}$. Finite dimensional vector spaces are replaced by 
modules defined by restriction of scalars from the one-dimensional free module, using the Frobenius endomorphisms of $\Z_{\rm max}$. The associated projective spaces are {\em finite} and  provide a mathematically consistent  interpretation of J.~Tits' original idea  of a geometry over the absolute point. The self-duality of the cyclic category and the cyclic descent number of permutations both acquire a geometric meaning.
\end{abstract}

\maketitle
\vspace*{2pt}
\tableofcontents  

\section{Introduction}

In this paper we establish a bridge between the combinatorial structure underlying cyclic homology and the $\lambda$-operations on one side and the framework of geometry in characteristic one on the other. 
The combinatorial system supporting cyclic homology and the $\lambda$-operations is best encoded by the cyclic category \cite{CoExt}  and its natural extension to the epicyclic category \cite{G1,bu1} which play an important   role in algebraic topology and algebraic $K$-theory (\cf\cite{good0}). In \cite{cycarch}, we showed the relevance of cyclic homology  of schemes and the $\lambda$-operations for the cohomological interpretation of the archimedean local factors of L-functions of arithmetic varieties, opening therefore the road to applications of cyclic homology in arithmetic. \newline
Mathematics in characteristic one has two algebraic incarnations: one is provided by the theory of semirings and semifields supporting tropical geometry and idempotent analysis, while the other one is centered on the more flexible notions of hyperrings and hyperfields on which certain number-theoretic constructions repose. In our  recent work \cite{CC1, CC2, CC3, CC4, CC5} we explained the relevance of these two algebraic theories to promote the development   of an absolute geometry.\newline
In this paper we provide the geometric meaning of the cyclic and the epicyclic categories in terms of a projective geometry in characteristic one and we supply the relation of the above categories with the absolute point. In \S \ref{sectF1} we show that the epicyclic category $\tilde \Lambda$ is isomorphic to a category $\cP_\F$ of projective spaces over the simplest {\em infinite} semistructure  of characteristic one, namely the semifield  $\F=\Z_{\rm max}:= (\Z\cup\{-\infty\},\text{max},+)$ of ``max-plus integers'' (here denoted multiplicatively). The objects of  $\cP_\F$   are  projective spaces $\P(E)$ where the \mods~$E$ over $\F$ are obtained by restriction of scalars from the one-dimensional free \mod~using the endomorphisms of $\F$. These endomorphisms form the multiplicative  semi-group $\N^\times$: for each integer $n\in\N^\times$ the corresponding endomorphism is the Frobenius $\ff_n$: $\ff_n(x):=x^n$ $\forall x\in \F$. Let denote by $\F^{(n)}$ the \mod~over $\F$ obtained from $\F$  by restriction of scalars using $\ff_n\in \End(\F)$, then for $n\geq 0$ the projective spaces $\P(\F^{(n+1)})$ provide the 
complete collection of objects of $\cP_\F$. 
The morphisms in $\cP_\F$ are  projective classes of {\em semilinear} maps $f$ of \mods~over $\F$ which fulfill the condition $f^{-1}(\{0\})=\{0\}$. One also derives the definition of a full (but not faithful) functor $\P: \cP_\F\longrightarrow \fin$ to the category of {\em finite sets} which associates to a \mod~$E$ over $\F$ the quotient space $\P(E)=(E\setminus \{0\})/\F^\times$ (\cf~Remark~\ref{rem} (a)). If one restricts the construction of the morphisms in $\cP_\F$ to maps which are {\em linear} rather than semilinear, one obtains a subcategory $\cP_\F^1\subset\cP_\F$ {\em canonically} isomorphic to the cyclic category $\Lambda$: the inclusion functor $\cP_\F^1\hookrightarrow \cP_\F$ corresponds to the inclusion of the categories $\Lambda\subset \tilde \Lambda$.\newline
 It is traditional to view the category of finite sets as the limit for $q=1$ of the category of finite dimensional vector spaces over a finite field $\F_q$ and the symmetric group $S_n$ as the limit case of the general linear group ${\rm GL}_n(\F_q)$. There is however one feature of the category of finite dimensional vector spaces over a field which is not preserved by this analogy, namely the self-duality
provided by transposition of linear maps. Indeed, the cardinality of the set of maps $\Hom_{\fin}(X,Y)$ between two finite sets is a highly asymmetric function of the sets, whereas for vector spaces over $\F_q$ the cardinality of $\Hom_{\F_q} (E_1,E_2)$  is the symmetric function $q^{n_1n_2}$, for $n_j=\dim E_j$ ($j=1,2$).\newline
The geometric interpretation provided in this paper 
of the epi/cyclic categories and of the functor $\P$ refines and clarifies the above correspondence.  In \S\ref{dual} we prove that the well known self-duality of the cyclic category is described by transposition of linear maps. On the other hand, the failure of the extension of the property of self-duality to the epicyclic category is explained by the fact that the transpose of a semilinear map fails to be semilinear when the associated morphism of fields is not surjective.  In our construction the semilinearity of the maps is encoded by the functor ${\rm Mod}: \cP_\F\longrightarrow \N^\times$ to the multiplicative mono\"{\i}d of natural numbers (viewed as a small category with a single object) which associates to a morphism $f$ in $\cP_\F$ the integer $n\in \N^\times$ such that $f(\lambda x)=\ff_n(\lambda)x$ $\forall\lambda \in \F$.  ${\rm Mod}$ also provides, using the functor $\P: \cP_\F\longrightarrow \fin$,   a geometric interpretation of the cyclic descent number of arbitrary permutations as the measure of their semilinearity: \cf~Proposition~\ref{propdescent}.\newline
One can finally formulate a mathematically consistent interpretation of  J.~Tits' original idea \cite{Tits} of a geometry over the absolute point which is provided in our construction by the data given by the category $\cP_\F$ ($\F = \Z_{max}$) and the functor $\P$. Notice that the cardinality of the set underlying the projective space $\P(\F^{(n+1)})$ is $n+1$ and that this integer coincides with the limit, for $q\to 1$, of the cardinality of the projective space $\P^n(\F_q)=\P(\F_q^{n+1})$. The fullness of the 
functor $\P$ shows in particular that any permutation $\sigma\in S_{n+1}$ arises from a geometric morphism of projective spaces over $\F$. \newline
Even though the above development of a (projective) geometry in characteristic one is formulated in terms of algebraic semistructures, in \S \ref{hypersec} we show how one can obtain its counterpart  in  the framework of hyperstructures  by applying a natural functor $-\otimes_{\B}\mathbf S$, where $\mathbf S$ is the smallest finite hyperfield of signs (\cf\cite{CC4}) that minimally contains the  smallest finite idempotent semifield $\B$. In \cite{CC5} we have shown that by implementing the theory of hyperrings and hyperfields  one can parallel successfully 
J.~M.~Fontaine's $p$-adic arithmetic theory of ``perfection'' and subsequent Witt extension by combining a process of dequantization (to characteristic one) and a consecutive Witt construction (to characteristic zero).  In view of the fact that this dequantization process needs the framework of hyperstructures to be meaningful, it seems evident that the arithmetical standpoint in characteristic one  requires a very flexible algebraic theory which encompasses semistructures. On the other hand, several successful developments of the theory of semirings in linear algebra and analysis show that the context of semistructures is already adequate for many  applications. The only reasonable conclusion one can draw is that for the general development of mathematics in characteristic one ought to keep both constructions available and select the most appropriate one in relation to the specific context in which each problem is formulated.

\section{The epicyclic category}

In this section we show that the notion of archimedean set and related category $\Arc$ (that we introduced in  \cite{topos}), provides a natural framework for  the definition of the variants $\Lambda_a$ (\cf \cite{bu1, good0}) of the cyclic category $\Lambda$  of \cite{CoExt} and of the  epicyclic category $\tilde \Lambda$ (originally due to Goodwillie). All these categories can be obtained by restricting to archimedean sets whose  underlying set is the set $\Z$ 
of integers with the usual total order. In this section we study the categories $\Arc$, $\Arc_a$ and $\Arc\triangleleft \N$ obtained by dropping the above restriction.

\subsection{The category $\Arc$ of Archimedean sets}

We recall from \cite{topos} the following notion
\begin{defn}\label{defnarc} An archimedean set is a pair $(X,\theta)$ of
a non-empty, totally ordered set $X$  and an order automorphism $\theta\in \Aut X$ such that $\theta(x)>x$, $\forall x\in X$ and fulfilling the following archimedean property:
\begin{equation*}
\forall x,y\in X, \ \exists n\in \N: \quad y\leq \theta^n(x).
\end{equation*}
\end{defn}

For each positive integer $a\in \N$ we introduce the following category $\Arc_a$
\begin{defn} \label{last} 
The objects of the category $\Arc_a$ are archimedean sets $(X,\theta)$;  the morphisms  
$f: (X,\theta)\to (X',\theta')$ in $\Arc_a$ are  equivalence classes of maps 
\begin{equation}\label{deff}
f:X\to X', \ \ f(x)\geq f(y) \  \  \forall x\geq y ;\qquad  f(\theta(x))=\theta'(f(x)),\quad \forall x \in X
\end{equation}
where the equivalence relation identifies two such maps $f$ and $g$ if there exists an integer $m\in \Z$ such that 
$g(x)=\theta'^{ma}(f(x))$, $\forall x\in X$.
\end{defn}

For $a=1$ we shall drop the index $1$: $\Arc_1=\Arc$ coincides with the category of archimedean sets.

\begin{prop}\label{propcorr0}
The full subcategory of $\Arc_a$ whose objects are the archimedean sets $(\Z,\theta)$, where $\Z$ is endowed with the usual order, is canonically isomorphic to the $a$-cyclic category $\Lambda_a$ considered in \cite{bu1, good0}. 
\end{prop}
\proof
One checks that the category of the archimedean sets such as $(\Z,\theta)$ is an extension (likewise $\Lambda_a$) of the small simplicial category $\Delta$ by means of a new generator $\tau_{n}$ of the cyclic group $C_{(n+1)a}=\Aut_{\Lambda_a}([n])$, for each $n\geq 0$, that fulfills the relations (\cf~ \cite{good0}, p. 235)
\begin{align*}
&&\tau_n^{(n+1)a}=id,&\notag\\
\tau_n\circ \sigma_0&=\sigma_n\circ \tau_{n+1}^2 &&&   \tau_n\circ \sigma_j &=\sigma_{j-1}\circ \tau_{n+1}, \qqq j\in \{1, \ldots, n\}\\
\tau_n\circ \delta_0 &=\delta_n,  &&&  \tau_n\circ \delta_j&=\delta_{j-1}\circ \tau_{n-1} \qqq j\in \{1, \ldots, n\}.
\end{align*}
\endproof
\begin{defn}\label{hatn} We denote with $\hat n=(\Z,\theta)$  the  archimedean set whose automorphism $\theta:\Z\to\Z$ is  defined by the translation $\theta(x)=x+n$, for a fixed  $n\geq 1$. 
\end{defn}
Such object $\hat n$ gives rise  to the object $[n-1]$ of  $\Lambda_a$:  the shifted indexing will be more convenient for our applications.

\subsection{The correspondences $\overline\Psi_k$}\label{sectcorr}

Let $(X,\theta)$ be an archimedean set and let $k>0$ be an integer. Then the pair $(X,\theta^k)$ is also an archimedean set that we denote as 
\begin{equation}
\label{psik}
\Psi_k(X,\theta):=(X,\theta^k).
\end{equation}
For $(X,\theta)$ and $(X',\theta')$ two archimedean sets and $f: (X,\theta)\to (X',\theta')$ a morphism in $\Arc$ connecting them,  one  has $f(\theta^k(x))=\theta'^k(f(x)),\  \forall x \in X$ ($k>0$ fixed). Thus $f$ defines a morphism $\Psi_k(f)\in \Hom_\Arc(\Psi_k(X,\theta),\Psi_k(X',\theta'))$.  However the two  maps $f$ and $\theta'\circ f$ which define the same morphism in the category $\Arc$ are in general no longer equivalent as morphisms $\Psi_k(X,\theta)\to\Psi_k(X',\theta')$. More precisely, one derives a correspondence
 $\overline\Psi_k: \Arc\DashedArrow[->,densely dashed] \Arc $ rather than a functor that satisfies the following properties

\begin{prop}\label{propcorr}
$(i)$~Let $h\in \Hom_\Arc((X,\theta),(X',\theta'))$, then for a fixed positive integer $k>0$ the set 
\begin{equation*}
\overline \Psi_k(h):=\{\Psi_k(f) \mid f\in h\}
\end{equation*}
is finite with exactly $k$ elements.\vspace{.1in}

$(ii)$~Let $g,h$ be composable morphisms in $\Arc$, then one has 
\begin{equation*}
 \overline\Psi_k(g\circ h):= \overline\Psi_k(g)\circ \overline\Psi_k(h)=\{u\circ v\mid u\in \overline\Psi_k(g) , \  v\in \overline\Psi_k(h) \}.
\end{equation*}
$(iii)$~For any positive integers $k,k'$:~ $\overline\Psi_k \circ \overline\Psi_{k'}=\overline\Psi_{kk'}$.

\end{prop}
\proof $(i)$~Let $f: (X,\theta)\to (X',\theta')$ be a morphism in $\Arc$, then the composite  $\theta'^k\circ f$ is equivalent to $f$ in the set  $\Hom_\Arc(\Psi_k(X,\theta),\Psi_k(X',\theta'))$, while the class of $f$ in $\Hom_\Arc((X,\theta),(X',\theta'))$ is represented by  $\theta'^m\circ f$, for $m \in \Z$.
Thus $ \overline\Psi_k(h)$ is the finite set of classes of  $\theta'^m\circ f$, for $m\in \{0, \ldots, k-1\}$. These elements are pairwise inequivalent since the maps $\theta'^m\circ f$ are pairwise distinct for $m\in \Z$.\newline
$(ii)$~Let $h\in \Hom_\Arc((X,\theta),(X',\theta'))$,  and $g\in \Hom_\Arc((X',\theta'),(X'',\theta''))$. Let  $g'\in g$ and $h'\in h$ be maps in the corresponding equivalence classes (fulfilling \eqref{deff}). Then $g'\circ h'\in g\circ h$, and one also has by construction 
 $\Psi_k(g'\circ h'):= \Psi_k(g')\circ \Psi_k(h') $.  By replacing $h'$ by $\theta'^n\circ h'$ and $g'$ by $\theta''^m\circ g'$ one substitutes $z=g'\circ h'$ with $\theta''^a\circ z$ with $a=n+m$ and only the class of $a$ modulo $k$ matters for the corresponding morphism from $\Psi_k(X,\theta)\to\Psi_k(X'',\theta'')$.\newline
$(iii)$~For any morphism $f: (X,\theta)\to (X',\theta')$  in $\Arc$ $f:X\to X'$ one easily check that $\Psi_k \circ \Psi_{k'}(f)=\Psi_{kk'}(f)$.\endproof

\subsection{Two functors $\Arc_a\longrightarrow\Arc_b$ when $b\vert a$}\label{sectdiv}

The correspondences $\overline\Psi_k: \Arc\DashedArrow[->,densely dashed] \Arc$ are best described in terms of two functors $P$ and $\Psi_k:\Arc_k\longrightarrow\Arc$, which we now describe in slightly more general terms.  

Let $a, b\in \N$: when $b\vert a$, the functor $P:\Arc_a\longrightarrow\Arc_b$ is the natural ``forgetful" functor. It is the identity on objects and associates to an equivalence class (Definition \ref{last}) of morphisms $f\in \Hom_{\Arc_a}((X,\theta),(X',\theta'))$ the unique class it defines in $ \Hom_{\Arc_b}((X,\theta),(X',\theta'))$.

The definition of $\Psi_k$ given in \S \ref{sectcorr} determines, for every positive integer $t>0$, a functor 
\begin{equation}\label{kmorpb}
   \Psi_k: \Arc_{kt} \longrightarrow \Arc_t
\end{equation}
this because the two  maps $f$ and $\theta'^{kt}\circ f$, which define the same morphism in the category $\Arc_{kt}$, are  equivalent as morphisms of the set $ \Hom_{\Arc_t}(\Psi_k(X,\theta),\Psi_k(X',\theta'))$. One thus obtains the following commutative diagram where the lower horizontal arrow is the correspondence  $\overline\Psi_k$
\begin{gather*}
 \,\hspace{130pt}\raisetag{-47pt} \xymatrix@C=35pt@R=35pt{
 \Arc_k\ar[dr]^-{\Psi_k} \ar[d]^{P}& \\
\Arc \ar@{-->}[r]^{\overline\Psi_k}  &  \Arc\\
}\hspace{140pt}
   \end{gather*}

\subsection{The category $\Arc\triangleleft \N$}

Proposition \ref{propcorr} describes an action by correspondences of the multiplicative mono\"{\i}d $\N^\times$ on $\Arc$. We shall now define the category $\Arc\triangleleft \N$ obtained as the semi-direct product of $\Arc$ by this action. Its objects are the same as those of $\Arc$. At the level of morphisms, instead, we adjoin, for each object $(X,\theta)$ of $\Arc$ and each positive integer $k$, a new morphism 
\begin{equation}\label{kmorpa}
   \psi_k: \Psi_k(X,\theta) \to (X,\theta)
\end{equation}
which fulfills the relations
\begin{equation}\label{kmorp1}
 f\circ   \psi_k= \psi_k\circ g \qqq g \in \Psi_k(f),\quad\forall f\in \Hom_\Arc((X,\theta),(X',\theta')) 
\end{equation}
and
\begin{equation}\label{kmorp2}
\psi_k\circ   \psi_{k'}= \psi_{kk'} \qqq k, k'>0.
\end{equation}
This construction is precisely achieved as follows

\begin{defn}\label{last2} The objects of the category $\Arc\triangleleft \N$ are the archimedean sets $(X,\theta)$;  the morphisms  
$f: (X,\theta)\to (X',\theta')$ in $\Arc\triangleleft \N$ are  equivalence classes of maps 
\begin{equation}\label{deff2}
f:X\to X', \ \ f(x)\geq f(y) \  \  \forall x\geq y ;\qquad \exists k\in \N, \  f(\theta(x))=\theta'^k(f(x)),\quad \forall x \in X
\end{equation}
where the equivalence relation identifies two such maps $f$ and $g$ if there exists an integer $m\in \Z$ such that 
$g(x)=\theta'^m(f(x))$, $\forall x\in X$.
\end{defn}
We check that the equivalence relation is compatible with the composition of maps in $\Arc\triangleleft \N$. 

Let  $f_j: (X,\theta)\to (X',\theta')$ ($j=1,2$) and $g: (X',\theta')\to (X'',\theta'')$ be two morphisms. One has $f_1\sim f_2~\Leftrightarrow~f_2(x)=\theta'^m(f_1(x))$, $\forall x\in X$ and for some $m\in\Z$. It follows (since $g$ fulfills \eqref{deff2}) that $g\circ f_2(x)=g(\theta'^m(f_1(x)))=\theta''^{km}g(f_1(x))$, 
thus $g\circ f_2\sim g\circ f_1$.\vspace{.05in}

Next proposition shows that the category $\Arc\triangleleft \N$ has exactly the expected properties of a semi-direct product of $\Arc$ by the correspondences $\bar\Psi_k$.

\begin{prop}\label{propsemi}
$(i)$~The category $\Arc$ is a subcategory of $\Arc\triangleleft \N$.\newline
$(ii)$~The map $\rho: \Hom_{\Arc\triangleleft \N}((X,\theta),(X',\theta')) \to \N$, $\rho(f)=k\in \N^\times$ such that \eqref{deff2} holds, describes at the morphisms level a functor 
\[
\rho: (\Arc\triangleleft \N)\longrightarrow\N^\times.
\]
$(iii)$~For any archimedean set $(X,\theta)$ and positive integer $k$, the identity map $id_X(x)=x$, $\forall x\in X$, defines a
morphism $ \psi_k\in \Hom_{\Arc\triangleleft \N}(\Psi_k(X,\theta),(X,\theta))$ which fulfills the relations \eqref{kmorp1}  and \eqref{kmorp2}.\newline
$(iv)$~Every morphism $f\in \Hom_{\Arc\triangleleft \N}((X,\theta),(X',\theta'))$ is  of the form 
\begin{equation*}
f=\psi_k\circ   h, \quad  k=\rho(f), \  \   h\in \Hom_\Arc(X, \Psi_k(X',\theta')).
\end{equation*}
\end{prop}
\proof $(i)$~The categories $\Arc$ and   $\Arc\triangleleft \N$ share the same objects and by construction one has $$\Hom_\Arc((X,\theta),(X',\theta'))\subset\Hom_{\Arc\triangleleft \N}((X,\theta),(X',\theta'))$$
$(ii)$~Since the action of $\theta'$ on $X'$ is free, the value of $k$ for which \eqref{deff2} holds is uniquely determined, moreover one checks easily that $\rho(f\circ g)=\rho(f)\rho(g)$.\newline
$(iii)$~The identity map $u=id_X$ fulfills $u\circ \theta^k=\theta^k\circ u$, thus defines a morphism 
$\psi_k$ in the set $\Hom_{\Arc\triangleleft \N}(\Psi_k(X,\theta),(X,\theta))$. By applying the definition of the equivalence relation for morphisms as in Definition \ref{last2}, one has $\psi_k\circ \theta^j\sim \psi_k$ for all $j\in \Z$. One thus obtains the equality $ f\circ   \psi_k= \psi_k\circ g$, $ \forall g \in \Psi_k(f)$. One checks easily the relations \eqref{kmorp1}  and \eqref{kmorp2} .\newline
$(iv)$~Let $f: (X,\theta)\to (X',\theta')$ in $\Arc\triangleleft \N$. Then $f$ determines an element $h\in \Hom_\Arc(X, \Psi_k(X',\theta'))$ whose definition depends upon the choice of $f$ in its class in $\Hom_{\Arc\triangleleft \N}((X,\theta),(X',\theta'))$. Replacing $f$ by $\theta'^j\circ f$ has the effect to replace $h$ by $\theta'^j\circ h$ whose class in $\Hom_\Arc(X, \Psi_k(X',\theta'))$ depends 
on the residue of $j$ modulo $k$.\endproof

\begin{prop}\label{propsemi1}
The full subcategory of $\Arc\triangleleft \N$ whose objects are the archimedean sets $(\Z,\theta)$ ($\Z$  endowed with the usual order) is {\em canonically isomorphic} to the epicyclic category $\tilde \Lambda$.
\end{prop} 
\proof By definition (\cf ~Definition~1.1 of \cite{bu1}), the epicyclic category $\tilde \Lambda$ is obtained by adjoining to the cyclic category $\Lambda$, new morphisms $\pi_n^k:[k(n+1)-1]\to [n]$ for $n\geq 0$, $k\geq 1$, which fulfill the following relations: 
\begin{enumerate}
\item $\pi_n^1={\rm id_n} $, $\pi_n^\ell \circ \pi_{\ell(n+1)-1}^{k}=\pi_n^{k\ell}$
\item $\alpha \pi_m^k=\pi_n^k {\rm Sd_k}(\alpha)$, for any $\alpha\in \Hom_\Delta([m],[n])$
\item $\tau_n \pi_n^k=\pi_n^k \tau_{k(n+1)-1}$.
\end{enumerate}
Here ${\rm Sd_k}:\Delta\longrightarrow \Delta$ is the barycentric subdivision functor which maps $[n-1]$ to $[kn-1]$ and a map $f$ to the $k$-fold concatenation ${\rm Sd_k}(f)=f\perp f \perp \ldots \perp f$. In terms of the archimedean sets $\hat n\sim [n-1]$ (\cf~Definition~\ref{hatn}) the map $f\in \Hom_\Delta([m-1],[n-1])$ lifts uniquely to a map $\tilde f: \Z\to\Z$ such that $\tilde f(x+m)=\tilde f(x)+n$ $\forall x\in \Z$. Notice that  $\tilde f$ agrees with $f $ on $\{0,1,\ldots,m-1\}$. Moreover $\tilde f$ is a morphism of archimedean sets and the class of
$\Psi_k(\tilde f)$ is the same as the class of  the $k$-fold concatenation $\widetilde{\rm Sd_k}(f)$. This shows that one obtains the required isomorphism of categories by extending the isomorphism of Proposition \ref{propcorr0} (for $a=1$) to the full subcategory of $\Arc\triangleleft \N$ by mapping the morphism $\psi_k\in  \Hom_{\Arc\triangleleft \N}(\Psi_k(\hat n),\hat n)$ to 
$\pi_{n-1}^k$.\endproof
\begin{rem}{\rm The epicyclic category $\tilde \Lambda$ used here is originally due to Goodwillie and described in  \cite{bu1}
but it does not correspond to the notion of  epicyclic space applied in \cite{good}.
}\end{rem}

\subsection{The functor $\cF:(\Arc\triangleleft \N)\longrightarrow\Se$}\label{sectF}

In the following we shall denote with $\mu_a$ the multiplicative group of $a$-th roots of unity in $\C$. By $\Se_a$ we denote the category of sets endowed with a free action of $\mu_a$, and with morphisms given by $\mu_a$-equivariant maps. For $(X,\theta)\in \text{Obj}(\Arc)$, we consider the orbit space of the action of $\theta^a$ on $X$:
\begin{equation*}
\cF_a(X,\theta):= X/\theta^{a\Z}
\end{equation*}
endowed with the free action of $\mu_a$ generated by the action of $\theta$ on $\cF_a(X,\theta)$.
\begin{prop}\label{propsemi2}
$(i)$~For $a=1$ one has a functor $\cF=\cF_1: (\Arc\triangleleft \N)\longrightarrow\Se$.\newline
$(ii)$~For any integer $a>1$, one has a functor $\cF_a:\Arc_a\longrightarrow\Se_a$.
\end{prop}
\proof $(i)$~Let $f:(X,\theta)\to (X',\theta')$ be a morphism in $\Arc\triangleleft \N$  (thus fulfilling \eqref{deff2}). Then given two points $x, y=\theta^m(x)$ on the same  orbit of the action of $\theta$ on $X$, the points $f(x)$ and $f(y)=\theta'^{km} (f(x))$ are on the same orbit of the action of $\theta'$ on $X'$. This shows that $\cF_1$ transforms a morphism $f\in \Hom_{\Arc\triangleleft \N}((X,\theta),(X',\theta'))$ into a map of sets, this association is also  compatible with composition of morphisms. \newline
$(ii)$~By definition of the equivalence relation as in  Definition \ref{last} for $f\in \Hom_{\Arc_a}((X,\theta),(X',\theta'))$, the induced map of sets $X/\theta^{a\Z}\to X'/\theta'^{a\Z}$ is independent of the choice of $f$ in its equivalence class. Moreover the equivariance condition  \eqref{deff} ensures that the induced map of sets  is $\mu_a$-equivariant. \endproof

We shall now follow the effect of the  functors $P$ and $\Psi_k$ as in \S \ref{sectdiv} in terms of the categories $\Se_a$. When $b\vert a$ ($a,b\in\N$), there is a canonical inclusion $\mu_b\subset \mu_a$. With $a=kb$ ($k\in\N$), the subgroup $\mu_b$ is the range of the group endomorphism $\mu_a\to\mu_a$ $u\mapsto u^k$. 
This determines a natural restriction functor 
\[
{\rm Res}: \Se_a\longrightarrow \Se_b
\]
 which does not alter the underlying set and restricts the action of the roots of unity $\mu_a$ to the subgroup $\mu_b$. This restriction functor corresponds to the  functor $\Psi_k:\Arc_a\longrightarrow\Arc_b$ of \eqref{kmorpb} \ie the following diagram commutes
\begin{gather*}
 \,\hspace{130pt}\raisetag{-67pt} \xymatrix@C=35pt@R=35pt{
\Arc_a\ar[d]_{\cF_a}\ar[r]^-{\Psi_k} &
\Arc_b\ar[d]^{\cF_b}& \\
\Se_a\ar[r]^-{{\rm Res}}  & \Se_b\\
}\hspace{100pt}
   \end{gather*}

To check the commutativity of the above diagram we note that the set underlying $\cF_b(\Psi_k(X,\theta))$ is the orbit space of the action of $\theta^{kb}$ on $X$
and this coincides with the set underlying $\cF_a(X,\theta)$.

Similarly, when $a=kb$, one has an ``extension of scalars" functor 
\[
-\times_{\mu_a} \mu_b: \Se_a\longrightarrow \Se_b
\]
 which associates to an object
$Y$ of $\Se_a$ its quotient $Y\times_{\mu_a} \mu_b$ for the action of the subgroup $\mu_k\subset \mu_a$. 
This functor corresponds to the  functor $P:\Arc_a\longrightarrow\Arc_b$ \ie the following diagram commutes:
\begin{gather*}
 \,\hspace{130pt}\raisetag{-47pt} \xymatrix@C=35pt@R=35pt{
\Arc_a\ar[d]_{\cF_a}\ar[r]^-{P} &
\Arc_b\ar[d]^{\cF_b}& \\
\Se_a\ar[r]^-{-\times_{\mu_a} \mu_b}  & \Se_b\\
}\hspace{100pt}
   \end{gather*}
\begin{rem}{\rm It is customary to interpret the category $\Se_a$ as the category of ``vector spaces over $\F_{1^a}$", where $\F_{1^a}$ plays the role of the limit for $q\to 1$ of the finite fields $\F_{q^a}$. However, this analogy has its limitations since for instance the classical duality between vector spaces over fields does not apply here since for ``vector spaces $V$ and $W$ over $\F_{1^a}$" of respective dimensions $n$ and $m$ the cardinality of the space of morphisms is
\begin{equation*}
  \# (\Hom_{\F_{1^a}}(V,W))=(am)^n
\end{equation*}
which is not a symmetric function of $n$ and $m$. In \S\ref{dual} we will explain how the classical duality is restored  for the cyclic category $\Lambda$ using the framework of projective geometry in characteristic one which we now describe. 
}\end{rem}

\section{$\tilde \Lambda$ and  projective geometry over the semifield $\Z_{max}$}\label{sectF1}

In ordinary projective geometry  the maps between projective spaces $\P(E_j)=(E_j\setminus \{0\})/K_j^\times$ ($j=1,2$) over fields $K_j$ are induced by {\em semilinear} maps of vector spaces $E_j$ (\cf\cite{Faure}).  Recall that a  map $f : E_1\to  E_2$ between two vector spaces is called semilinear if it is additive and if there exists a homomorphism of fields $\sigma: K_1\to K_2$ such that $f(\lambda x)=\sigma(\lambda)f(x)$ $\forall \lambda \in K_1$ and $\forall x\in E_1$. This notion extends verbatim to the context of semifields where by a semifield we mean a commutative semiring $K$ in which the non-zero elements form a group under multiplication (\cf \cite{Golan}, 4.25, p. 52) and in a semiring the existence of an additive inverse is no longer required (\cf \cite{Golan}, I).

By a \mod~$E$ over a semifield $K$ we mean (\cf\cite{Golan}, Chapter 14, \cf \cite{Gondran}, Chapter 5) a commutative mono\"{\i}d $(E,+)$ with additive identity $0\in E$, endowed with an action of $K$ such that 
$\forall\lambda,\mu\in K$ and $\forall x,y\in E$, one has 
 \begin{equation}\label{vect0}
\lambda (x+y)=\lambda x+\lambda y, \  \  (\lambda+\mu) x=\lambda x+ \mu x, \  \  (\lambda\mu) x=\lambda(\mu x), \  \  0 \ x=0, \  \  1\ x=x.
\end{equation}
A  map $f : E_1\to  E_2$ between two \mods~over semifields $K_j$ is called semilinear if it is additive and if there exists a homomorphism of semifields $\sigma: K_1\to K_2$ such that $f(\lambda x)=\sigma(\lambda)f(x)$ $\forall \lambda \in K_1$ and $\forall x\in E_1$. Two such maps $f,f'$ are projectively equivalent when there exists $\lambda\in K_2$, $\lambda \neq 0$, such that $f'(x)=\lambda f(x)$,  $\forall x\in E_1$.
\begin{defn} \label{defnpro} Let $\proj$ be the category whose objects are pairs $(K,E)$ made by a semifield $K$ and a \mod~$E$ over $K$ and whose morphisms $(K_1, E_1)\to (K_2, E_2)$ are pairs $(\sigma, h)$ where $\sigma: K_1\to K_2$ is a semifield homomorphism and $h$ is a projective class of additive maps $f:E_1\to E_2$ such that such that $f^{-1}(\{0\}) = \{0\}$ and  $f(\lambda x)=\sigma(\lambda)f(x)$ $\forall \lambda \in K_1$, $\forall x\in E_1$.
\end{defn}

\begin{prop}\label{defnP} 
The assignment  which maps a pair $(K,E)\in\text{Obj}(\proj)$  to the corresponding projective space (as set) $\P(E):=(E\setminus \{0\})/K^\times$ and a morphism in $\proj$ to the induced map of sets defines a covariant functor $\P: \proj\longrightarrow \Se$. 
\end{prop}
\proof Let $f : E_1\to  E_2$ be a semilinear map such that $f^{-1}(\{0\}) = \{0\}$. For $x\in E_1\setminus \{0\}$ the class of $f(x)\in (E_2\setminus \{0\})/K_2^\times$ does not change if one replaces $x$ by $\lambda x$ for $\lambda \in K_1^\times$ or if one replaces $f$ by $\mu f$ for $\mu \in K_2^\times$.
\endproof
The assignment which associates to an object $(K,E)$ of $\proj$ the semifield $K$ defines a functor 
${\rm Mod}$  from  $\proj$ to the category of semifields.
\begin{defn} Let $K$ be a semifield. We denote by $\P_K$  the full subcategory of $\proj$ whose objects are   \mods~over $K$. We denote by  $\P_K^1$ the subcategory of $\P_K$  with the same objects as $\P_K$ and whose morphisms are pairs $(\sigma,f)$ where $\sigma$ is the identity on $K$ (\ie morphisms in $\P_K^1$ are given by projective classes of linear maps). 
\end{defn}
Thus by definition the objects of $\P_K$ are the objects $X=(K,E)$ of $\proj$ such that ${\rm Mod}(X)=K$ and the morphisms $\alpha$ of $\P_K^1$ are such that ${\rm Mod}(\alpha)={\rm id}_K$.

Next, we recall the most used definition of rank of a \mod~(\cf~\eg \cite{Golan}, Chapter 14, page 153) and we also introduce the notion of free rank that generalizes, in the context of semistructures, the classical notion of the largest cardinality of a ``free system''.
\begin{defn} Let $E$ be a finitely generated  \mod~over a semifield $K$\newline
$a)$~the rank $\rk(E)$ is the smallest integer $n$ such that there exists a set of generators of $E$ of cardinality $n$.\newline
$b)$~the free rank $\ruk(E)$ is the largest integer $n$ such that there exists a free \modsub~of $E$ of rank $n$.
\end{defn}
Following Exercise 16 A II.181 of \cite{bbk} one can show that if $\rk(E)$ and $\ruk(E)$ are  both finite then $\ruk(E)\leq\rk(E)$.

\subsection{The semifield $\B$ and the simplicial category $\Delta$}

Unlike the classical case of vector-spaces over fields, a finitely generated \mod~$E$ over a semifield $K$ is not necessarily isomorphic to $K^n$, for some $n$. This change of behavior arises already in the simplest example of the idempotent semifield $K=\B=(\{0,1\},+,\cdot)$.
Here the term ``idempotent" means that $x+x=x$, $\forall x\in K$, or equivalently that $1+1=1$ which we view as reflecting the fact that one works in characteristic one. It is known (\cf~\cite{Golan}, Chapter 4, 4.28) that $\B$ is the only finite idempotent semifield.

 Let $E$ be a \mod~over $\B$. Since $1+1=1$ in $\B$ it follows from \eqref{vect0} that $x+x=x$, $\forall x\in E$ so that, as a mono\"{\i}d, $E$ is idempotent and we shall use the notation $x+y=x\vee y$ for the sum of two elements of $E$. The canonical preorder of the commutative mono\"{\i}d $E$, defined by (\cf \cite{Gondran} 3.3, page 12)
 \begin{equation*}
x\leq y \iff \exists z \in E, \ y=x\vee z
\end{equation*}
is then an order relation (\cf \cite{Gondran}, Proposition 3.4.5) and moreover one has $x\leq y\iff x\vee y=y$.

The  partially ordered set $E$ is a semi-lattice with a smallest element $0$ and the join of any two elements $x,y\in E$ is  $x\vee y$.
Conversely, given a semi-lattice  $X$, one defines a \mod~$X^ \vee$ over $\B$ by adjoining to $X$ a smallest element as follows 
\begin{defn} Let $X$ be a semi-lattice, one lets $X^ \vee=X\cup \{0\}$ be the set endowed with the following binary operation $\vee$ 
 \begin{equation}\label{vect0.5}
x\vee y :={\rm join}(x,y) \qqq x,y \in X, \qquad 0 \vee x=x\vee 0=x \qqq x\in X^ \vee.
\end{equation}
\end{defn}
It is easy to see that the last two equations of \eqref{vect0} uniquely define the action of $\B$ on $X^ \vee$ and that this action fulfills the other equations of \eqref{vect0} since $x\vee x =x$ $\forall x\in X^ \vee$.

\begin{defn} \label{select}(\cf\cite{Gondran} p. 5) 
A commutative mono\"{\i}d $E$ is {\em selective} if and only if one has $x+y\in \{x,y\}$, $\forall x,y\in E$.
\end{defn}
By \cite{Gondran}, Proposition 3.4.7, a commutative mono\"{\i}d is selective if and only if it is idempotent and the canonical order is total.

The following statement determines a  complete list of the finitely generated \mods~of free-rank one over $\B$ and their categorical interpretation 
\begin{prop} \label{propfin} $(i)$~Let $E$ be a $\B$-\mod, then $E$ is selective if and only if $\ruk(E)=1$.

$(ii)$~There exists a unique, up to canonical isomorphism, $\B$-\mod~$E=\B^{(n,1)}$ such that $\rk(E)=n$ and 
$\ruk(E)=1$. One has $\B^{(n,1)}=X^ \vee$, for $X=\{1,\ldots ,n\}$ as a totally ordered set.

$(iii)$~The following properties hold for the  \mods~$\B^{(n,1)}$:

$a)$~$\B^{(n,1)}$ is of minimal cardinality among the $\B$-\mods~of rank $n$.  \newline 
$b)$~$\B^{(n,1)}$ is a projective 
\mod~over $\B$ equal to the range of a projection matrix $P\in M_n(\B)$, $P^2=P$. \newline
$c)$~The semiring of endomorphisms $\End_\B(\B^{(n,1)})$ is isomorphic to $PM_n(\B)P$.

$(iv)$~The simplicial category $\Delta$ is canonically isomorphic to the full subcategory $\cP_\B\subset\P_\B$  whose objects are the  \mods~$\B^{(n,1)}$ and the morphisms are (projective classes of) linear maps $f$ such that $f^{-1}(\{0\})=\{0\}$.
\end{prop}

\proof $(i)$~Let $E$ be a $\B$-\mod~such that $\ruk(E)=1$, then $x\vee y \in \{x,y\}$ for any two non-zero distinct elements $x,y\in E$ since otherwise one could construct a free \modsub~of $E$ of rank two.  It follows that $E$ is selective. Conversely if $E$ is selective it does not contain a copy of $\B^2$ and thus $\ruk(E)=1$.

$(ii)$~The rank of $E=X^ \vee$ for $X$ totally ordered is $\rk(E)={\rm card}X$. Thus if $\rk(E)=n$ one has $E=\B^{(n,1)}=X^ \vee$ with $X=\{1,\ldots ,n\}$ as  a totally ordered set.\newline
$(ii)$~Let us show that $(a)$, $(b)$, $(c)$ hold. \newline
$a)$~The cardinality of $\B^{(n,1)}$ is $n+1$ and is the minimal cardinality among $\B$-\mods~of rank $n$, since any such  \mod~contains at least $0$ and the $n$ generators. Let $E$ be a $\B$-\mod~of rank $n$, and assume that the cardinality of $E$ is $n+1$, then $x\vee y \in \{x,y\}$ for any $x,y\in E$, since otherwise one could remove $x\vee y$ from the set of generators. Thus it follows that $E$ is selective and thus $E=\B^{(n,1)}$.

$b)$~One has 
 \begin{equation}\label{projP}
\B^{(n,1)}= P(\B^n), \qquad P=\begin{pmatrix}
   1   &   0 &0& \ldots &   0 \\
    1  &  1&0& \ldots &   0 \\
  1  &  1&1& \ldots &   0 \\
& \ldots && \ldots  &    \ldots \\
1&1&1& \ldots &   1 
\end{pmatrix}\in M_n(\B), \qquad P^2=P
\end{equation}
since the image $P(\{e_j\})$ of the canonical basis $\{e_j\}$ of $\B^n$ is the decreasing sequence
$
P(\{e_j\})=\{\vee_{i\leq j} e_i\}
$. This shows that $\B^{(n,1)}$ is a finitely generated and projective  \mod~of rank $\leq n$. By \eqref{vect0.5} any subset of $X=\{1,\ldots ,n\}$ is stable under the binary operation $\vee$ and this shows that the rank of  $X^ \vee$ is equal to $n$.\newline
$c)$~The set of endomorphisms $\End_\B(\B^n)$  forms a semiring isomorphic to the semiring of matrices $M_n(\B)$  (\cf\cite{Golan}). Given $T\in\End_\B(\B^{(n,1)})$ the composite $T\circ P$ defines an element of $\End_\B(\B^n)$.

$(iii)$~There is a natural functor which associates to   a totally ordered set $X$ the \mod~$X^\vee$ over $\B$. The morphisms 
$f\in\Hom_\B( X^\vee, Y^\vee)$ such that $f^{-1}(\{  0\})=\{  0\}$ are  the non-decreasing maps from $X$ to $Y$. 
By applying this construction to  the skeleton  of the category of totally ordered finite sets we obtain a functor from the simplicial category $\Delta$ to the  category $\cP_\B$. This  functor is fully faithful and hence it defines an  isomorphism of categories $\Delta\sim \cP_\B$. \endproof
\begin{rem}
{\rm $a)$~Claim $(iv)$ of Proposition \ref{propfin} does not change if one replaces linear maps by projective classes of linear maps since $\B^\times =\{1\}$. In fact, since $\B$ has no non-trivial endomorphism, linearity and semilinearity are equivalent notions in this context.

$b)$~The \modsubs~of rank $k$ of $X^ \vee$ are determined by the subsets of $X$ of cardinality $k$. Thus the cardinality of the set of \modsubs~of rank $k$ of $\B^{(n,1)}$ is given by the binomial coefficient ${n}\choose{k}$. This is in agreement with the limit, as $q\to 1$, of the cardinality of the Grassmannian of vector subspaces of dimension $k$ in 
an $n$-dimensional vector space over the finite field $\F_q$.

$c)$~The free \mod~$\B^n$ has cardinality $2^n$ which grows exponentially with $n$ while the cardinality of $\B^{(n,1)}$  is $n+1$ which is linear in $n$. 
}\end{rem}

\subsection{\Mods~over $\Z_{max}$ and archimedean sets}

Let $\F=\Z_{\rm max}:= (\Z\cup\{-\infty\},\text{max},+)$ be the  semifield of tropical integers: we shall denote it multiplicatively, thus the elements of $\F$ are either $0$ or a power $u^n$ for $n\in \Z$. The idempotent addition $\vee$ is such that $u^n\vee u^m=u^k$, with $k=\sup( n,m)$. The multiplication is the usual one: $u^n u^m=u^{n+m}$. $\F$ is isomorphic to the 
sub-semifield of $\rmax$ generated by an element $>1$ of  $\rmax$.  

In this section we interpret the category $\Arc$ and the functor $\cF=\cF_1:\Arc \longrightarrow \Se$ in terms of the category $\P_\F^1$ of \mods~over $\F=\Z_{\rm max}.$ \newline
An archimedean set $(X,\theta)$ defines a \mod~$(X,\theta)^ \vee$ over $\F$ as follows:  
\begin{prop}\label{my check} Let $(X,\theta)$ be an archimedean set. Let $(X,\theta)^ \vee=(X\cup \{0\},\theta)$ be endowed with the binary operation 
 \begin{equation}\label{vect1}
x\vee y :=\sup(x,y) \quad \forall x,y \in X, \qquad 0 \vee x=x\vee 0=x\quad \forall x\in X\cup\{0\}.
\end{equation}
The action of $\F$ on $(X,\theta)^ \vee$  given by   
 \begin{equation}\label{vect2}
u^n x :=\theta^n(x) \quad \forall x \in X, \ n \in \Z,\qquad 0 x=0 \quad \forall x\in X\cup\{0\}.
\end{equation}
endows $(X,\theta)^ \vee$ with the structure of \mod~over $\F$.
\end{prop}
\proof
The  condition $\theta(x)>x$, $\forall x\in X$, of Definition \ref{defnarc} shows that 
 \begin{equation*}
(a\vee b) x =a x\vee b x \qqq x \in X\cup\{0\}, \  a, b \in \F.
\end{equation*}
Moreover the linearity property
 \begin{equation*}
a(x\vee y)  =a x\vee a y \qqq x, y \in X\cup\{0\}, \  a \in \F
\end{equation*}
holds since $\theta$ is an order automorphism.\endproof
\begin{prop} \label{natfun}
There is a fully faithful functor $\vee: \Arc \longrightarrow \P_\F^1$, $(X,\theta)\mapsto (X,\theta)^ \vee$,  mapping morphisms in $\Arc$ to projective classes of linear maps $f$ such that $f^{-1}(\{0\})=\{0\}$.
\end{prop}
\proof We know already that if $(X,\theta)$ is an archimedean set then $(X,\theta)^ \vee$ is a \mod~over $\F$. Next we show how  to define the functor $f\mapsto f^\vee$ on morphisms. Let 
$f:(X,\theta)\to (X',\theta')$ be a morphism in $\Arc$ (thus fulfilling \eqref{deff}), then we extend $f$ at $0$ by  $f(0)=0$ and obtain an $\F$-linear map $f^ \vee: (X,\theta)^ \vee\to (X',\theta')^ \vee$.  By construction, one has $f^{-1}(\{0\})=\{0\}$. Moreover, replacing $f$ by $\theta'^m\circ f$ does not alter the projective class of $f^ \vee: (X,\theta)^ \vee\to (X',\theta')^ \vee$ since it replaces $f^ \vee$ by $u^m f^ \vee$. In this way one obtains a functor $\Arc\longrightarrow \P_\F^1$. It is faithful by construction, so it remains to show that it is full.  Let $h: (X,\theta)^ \vee\to (X',\theta')^ \vee$ be an $\F$-linear map such that $h^{-1}(\{0\})=\{0\}$. Then the restriction of $h$ to $X$ defines a map $f:(X,\theta)\to (X',\theta')$ that fulfills \eqref{deff}. This shows the required surjectivity of the functor on morphisms. \endproof

\begin{rem} {\rm
Any non-zero morphism $h: (X,\theta)^ \vee\to (X',\theta')^ \vee$ fulfills $h^{-1}(\{0\})=\{0\}$. Indeed, assume that $h^{-1}(\{0\})$ contains an element $x\in X$, then we prove that $h(y)=0$ $\forall y\in X$. The archimedean property shows that for any $y\in X$ there exists an integer $n$ such that $y\leq \theta^n(x)$. It follows that $h(y)\leq \theta'^n(h(x))=0$ and thus $h(y)=0$.
}\end{rem}

\subsection{Geometric interpretation of the functor $\cF_1$}
 In the above geometric terms, the functor $\cF_1:\Arc\longrightarrow \Se$ is a special case of the  functor $\P:\proj\longrightarrow \Se$ of Proposition \ref{defnP}, \ie
the following diagram commutes:
 \begin{gather*}
 \,\hspace{130pt}\raisetag{-47pt} \xymatrix@C=35pt@R=35pt{
\Arc\ar[dr]_{\cF_1}\ar[r]^-{\vee} &
 \P^1_\F\ar[d]^{\P}& \\
  & \Se \\
}\hspace{100pt}
   \end{gather*}

One may wonder what geometric structure remains after passing from a \mod~$E$ to the set $\P(E)$. In ordinary projective geometry (where $E$ would be a vector space over a field) this structure is given by the map $(x,y)\mapsto  \ell(x,y)$ which associates to a pair of distinct points of a projective space the line  determined by them. Then the axioms of projective geometry characterize the obtained structures in the Desarguesian case.\newline
In the above framework of projective geometry over $\F=\Z_{max}$, there is a similar geometric structure: the 
``abstract circle" of \cite{Moerdjik}.  By definition an
abstract circle $C$ is given by the following data
\begin{equation}
\label{ abstrcirc}
C=(P,S,\partial_0,\partial_1,0,1,*,\cup)
\end{equation}
here $P$ and $S$ are sets, $\partial_j:S\to P$ are maps as well as $P\ni x\to 0_x\in S$ and $P\ni x\to 1_x\in S$, $*:S\to S$ is an involution, and $\cup$ is a partially defined map from a subset of $S\times S$ to $S$. Here $P$ plays the role of the set of points of the geometry while $S$ plays the role of the set of lines, or rather ``segments". In order to qualify as an abstract circle the data \eqref{ abstrcirc} have to fulfill certain axioms (\cf\cite{Moerdjik, topos}). It follows from these axioms  that given two points $x\neq y\in P$ there exists a unique segment $s\in S$ such that $\partial_0 (s)=x$, ~$\partial_1 (s)=y$. 

 As shown in \cite{topos} there is a natural functor $\frak Q$ which associates to an object $(X,\theta)$ of the category $\Arc$  an abstract circle $X/\theta$
and establishes in this way  an equivalence of categories.

The abstract circle $X/\theta=(P,S,\partial_0,\partial_1,0,1,*,\cup)$ associated to an archimedean set $(X,\theta)$ is obtained as follows: 

$-$~$P:=X/\sim$~is the orbit space for the action of $\Z$ on $X$ given by powers of $\theta$, \ie  $P=\P(E):=(E\setminus\{0\})/\F^\times$ for $E$ \mod~over $\F$.\newline
$-$~$S$ is  the orbit space for the action of $\Z$ on the set of pairs $(x,y)\in X^2$, with $x\leq y\leq \theta(x)$.\newline
$-$~$\partial_0 (x,y)=x$, ~$\partial_1 (x,y)=y$.\newline
$-$~$0_x=(x,x)$,~$1_x=(x,\theta(x))$.\newline
$-$~$(x,y)^*=(y,\theta(x))$.\newline
$-$~$(x,y)\cup (y,z)=(x,z)$ provided that $x\leq y \leq z\leq \theta(x)$.
\vspace{.05in}

It remains to be seen how to relax the conditions fulfilled by the morphisms in the category of abstract circles so that the above discussion extends to the category $\Arc\triangleleft \N$ and the diagram below commutes

 \begin{gather*}
 \,\hspace{130pt}\raisetag{-47pt} \xymatrix@C=35pt@R=35pt{
\Arc\triangleleft \N\ar[dr]_{\cF}\ar[r]^-{\vee} &
 \P_\F\ar[d]^{\P}& \\
  & \Se \\
}\hspace{100pt}
   \end{gather*}

\subsection{Restriction  of scalars}

In this section we implement the semigroup of endomorphisms of the semifield $\F=\Z_{max}$ to define a functor restriction  of scalars for \mods~$E$ over $\F$. To each integer $n\in \N^\times$ corresponds an endomorphism $\ff_n\in \End(\F)$ given by $\ff_n(x):=x^n$ $\forall x\in\F$. Moreover one has
\begin{lem}\label{lemfrob}
The map $ \N^\times\to \End(\F)$, $ n\mapsto \ff_n$ is an isomorphism of semigroups.
\end{lem} 
Let $E$ be a \mod~over $\F$, and $n\in \N^\times$.
Since $\ff_n:\F\to \F$ is a homomorphism of semifields, one can associate to $E$ the \mod~$F_n(E)$ over $\F$ with the same underlying additive structure but with a 
re-defined multiplication by elements of $\F$ as follows
$$
\lambda . \xi := \ff_n(\lambda) \xi.
$$
Since $\ff_n$ is not surjective, this restriction of scalars fails to pass unambiguously to projective classes of $\F$-linear maps $f$ with $f^{-1}(\{0\})=\{0\}$. Indeed, the ambiguity is retained by the group $\F^\times/\ff_n(\F^\times)\sim \mu_n$.

\subsection{Semilinear maps and $\Arc\triangleleft \N$}

Next, we extend Proposition \ref{natfun} to the category $\Arc\triangleleft \N$.
\begin{prop} \label{natfun1} The functor $\vee: \Arc\triangleleft \N \longrightarrow \P_\F$, $(X,\theta)\mapsto (X,\theta)^ \vee$ is fully faithful.
\end{prop}
\proof Recall (\cf~ Definition \ref{last2}) that the objects of the category $\Arc\triangleleft \N$ are archimedean sets while the morphisms 
are equivalence classes of maps which fulfill \eqref{deff2}. The condition $ f(\theta(x))=\theta'^k(f(x))$, $ \forall x \in X$ implies that extending $f$ by $f(0)=0$ one obtains an $\F$-semilinear map $f^ \vee: (X,\theta)^ \vee\to (X',\theta')^ \vee$,
with $ f^ \vee(\lambda \xi)=\ff_k(\lambda)f^ \vee(\xi)$ $\forall \lambda\in \F$ and $\xi\in (X,\theta)^ \vee$. Since  any morphism $\sigma:\F\to \F$ is an $\ff_k$ for some $k\in \N^\times$ (\cf~Lemma \ref{lemfrob}), the proof  of Proposition \ref{natfun} applies verbatim to show that the obtained functor is full and faithful from the category $\Arc\triangleleft \N$ to the category of 
\mods~over $\F$ with morphisms given by projective classes of semilinear maps $f$ such that $f^{-1}(\{0\})=\{0\}$.
\endproof

\subsection{The epicyclic category $\tilde\Lambda$ and projective geometry over $\Z_{max}$} 

We first investigate the structure of the \mods~$E$ over $\F$ obtained from the one dimensional free \mod~by restriction of scalars using the endomorphisms of $\F$.

\begin{prop}\label{propepi}
$(i)$~The \mods~over $\F$ of the form $(\Z,\theta)^\vee$ (\ie coming from archimedean sets of the form $(\Z,\theta)$) are obtained by restriction of scalars and they are of the form $F_n(\F)=\F^{(n)}$ where the integer $n\in \N$ is such that $\theta(x)=x+n$ $\forall x\in \Z$. 

$(ii)$~The functor $E\longrightarrow \P(E)$ establishes a bijection between the \modsubs~of $\F^{(n)}$ and the subsets of $\P(\F^{(n)})$.

$(iii)$~Let $E$ be a \modsub~of $\F^{(n)}$, then $E$ is isomorphic to $\F^{(k)}$ where $k$ is the rank of $E$.

\end{prop}
\proof $(i)$~Follows from the definition of the restriction of scalars. 

$(ii)$~One has by construction $\P(E)\subset \P(\F^{(n)})$ and this gives an injection  between the \modsubs~of $\F^{(n)}$ and the subsets of $\P(\F^{(n)})$. To prove that this map is surjective it is enough to show that given a subset $Y\subset \P(\F^{(n)})$ with $k$ elements there exists a morphism $f: \F^{(k)}\to \F^{(n)}$ such that the range of $\P(f)$ is $Y$. This statement will also prove $(iii)$ provided that $f$ is an isomorphism with its range. From $(i)$ one has $\F^{(n)}=(\Z,\theta)^\vee$ where $\theta(x)=x+n$ $\forall x\in \Z$. The subset $\{0,\ldots, n-1\}\subset \Z$ is a fundamental domain for the action of $\theta$ and gives an identification $\P(\F^{(n)})\sim \{0,\ldots, n-1\}$. One has $Y=\{y_0,\ldots,y_{k-1}\}\subset \{0,\ldots, n-1\}$, where $y_0<y_1<\cdots < y_{k-1}$. Let
$$
f(x):=y_{\bar x}+n E(x/k) \qqq x\in \Z
$$
where $\bar x\in \{0,\ldots, n-1\}$ is the residue of $x$ modulo $n$ and $E(z) $ is the integer part of $z\in \R$. Then $f:\Z\to \Z$ is an increasing map and fulfills $f(x+k)=f(x)+n$ $\forall x\in\Z$. Thus $f$ defines an injective morphism $f: \F^{(k)}\to \F^{(n)}$ such that the range of $\P(f)$ is $Y$.\endproof

Next statement shows that the epicyclic category  $\tilde\Lambda$ encodes  projective geometry over the semifield $\F=\Z_{max}$, where the projective spaces $\P(E)$  are constructed from the \mods~$E=\F^{(n)}$ over $\F$.
\begin{thm}\label{thmepi}
$(i)$~The epicyclic category $\tilde\Lambda$ is canonically isomorphic to the full subcategory $\cP_\F\subset\P_\F$ whose objects are obtained from the one dimensional free \mod~$\F=\Z_{max}$ by restriction of scalars using the endomorphisms of $\F$.

$(ii)$~The cyclic category $\Lambda\subset \tilde \Lambda$ is isomorphic to the subcategory $\cP^1_\F\subset\cP_\F$  with the same objects of $\cP_\F$  and whose morphisms are induced by  linear maps.  
\end{thm}
\proof $(i)$~The statement follows from Proposition \ref{natfun1}, Proposition \ref{propsemi1} and Proposition \ref{propepi}. Thus the object $[n-1]$ of $\tilde\Lambda$ corresponds canonically to $F_n(\F)=\F^{(n)}$ and Proposition \ref{natfun1} determines the canonical isomorphism. \newline
$(ii)$~follows from Proposition \ref{natfun}.
\endproof

Next, we investigate how the inclusion  $\Delta\subset \Lambda$ of the  simplicial category into the cyclic category arises from extension of scalars from $\B$ to $\F=\Z_{max}$. First, we need to relate the $\F$-\mod~$\F^{(n)}$ with the \mod~$\B^{(n,1)}$ ($n\in\N$).  Both pairs $(\B,\B^{(n,1)})$ and 
$(\F,\F^{(n)})$ are objects of the category $\proj$ as in  Definition \ref{defnpro}. 

Let $\iota:\B\to \F$ be the unique homomorphism of semifields. By construction $\F^{(n)}=(\Z,\theta)^\vee$ where  $\theta(x)=x+n$ $\forall x\in \Z$. 

Let  $\iota_n:\B^{(n,1)}\to \F^{(n)}$ be $\gamma^\vee$ where $\gamma$ is the unique increasing map which identifies the finite ordered set $\B^{(n,1)}\setminus \{0\}$ with the subset $\{0, \ldots,n-1\}\subset \Z$.

 \begin{prop}
$(i)$~The pair $(\iota,\iota_n)$ defines a morphism in $\proj$.

$(ii)$~Let $f\in\Hom_{\cP_\B}(\B^{(n,1)},\B^{(m,1)})$. Then there exists a unique $\tilde  f\in \Hom_{\cP^1_\F}(\F^{(n)},\F^{(m)})$ such that the following diagram commutes in $\proj$:
\begin{gather}
\label{funpp}
 \,\hspace{130pt}\raisetag{-27pt} \xymatrix@C=35pt@R=35pt{
\F^{(n)}\ar[r]^-{\tilde  f}  & \F^{(m)}\\
\B^{(n,1)}\ar[u]_{(\iota,\iota_n)}\ar[r]^-{f} &
\B^{(m,1)}\ar[u]^{(\iota,\iota_m)}& \\
}\hspace{100pt}
   \end{gather}

$(iii)$~The functor $\cP_\B\longrightarrow \cP^1_\F$, $f\mapsto \tilde f$,  corresponds to the canonical inclusion  $\Delta\subset \Lambda$. 
\end{prop}
\proof $(i)$~By construction $\iota_n:\B^{(n,1)}\to \F^{(n)}$ is additive and $\B$-linear, thus the pair $(\iota,\iota_n)$ defines a morphism
$
(\iota,\iota_n)\in \Hom_\proj((\B,\B^{(n,1)}),(\F,\F^{(n)}))
$
which is also represented for each $k\in \Z$ by the projectively equivalent pair $(\iota,u^k\iota_n)$.

$(ii)$~We identify $\B^{(n,1)}\setminus \{0\}$ with the subset $\{0, \ldots,n-1\}\subset \Z$. The morphism $f\in\Hom_{\cP_\B}(\B^{(n,1)}, \B^{(m,1)})$ is given by a unique non-decreasing map (we still denote it by $f$)
$f:\{0, \ldots,n-1\}\to \{0, \ldots,m-1\}$. To prove the existence, one defines the map $g:\Z\to \Z$ by
 \begin{equation}\label{gmap}
g(j+kn)=f(j)+km \qqq j\in \{0, \ldots,n-1\}, \ k \in \Z.
\end{equation}
One then gets that $g^\vee\in \Hom_{\cP^1_\F}(\F^{(n)},\F^{(m)})$ and that the diagram \eqref{funpp} is commutative which proves the existence of $g=\tilde f$. To prove the uniqueness of $\tilde f$ note that every non-zero element of $\F^{(n)}$ is of the form $x=u^k\iota_n(y)$ for some $y\in \B^{(n,1)}$. Thus if the diagram \eqref{funpp}, with $h$ instead of $\tilde f$, commutes in $\proj$, there exists $\ell \in \Z$ such that 
$
h(j+kn)=f(j)+km+\ell m \quad \forall j\in \{0, \ldots,n-1\}, \ k \in \Z
$
and this shows that $h$ is in the same projective class as the above $g$.

$(iii)$~By construction  $g=\tilde f$ fulfills \eqref{gmap} 
and this corresponds to the canonical embedding  $\Delta\subset \Lambda$. \endproof

\begin{rem}\label{rem}{\rm $a)$~When applied to the morphism $(\iota,\iota_n)\in \Hom_\proj((\B,\B^{(n,1)}),(\F,\F^{(n)}))$ the functor 
$\P$ of Proposition \ref{defnP} determines a {\em bijection} 
$$
\P(\B^{(n,1)})\stackrel{\P(\iota,\iota_n)}{\to} \P(\F^{(n)}).
$$
Thus, as a set, the projective space does not change by implementing an extension of scalars from $\B$ to $\F$ and moreover it remains {\em finite} of cardinality $n$. One derives the definition of a full functor
\[
\P: \cP_\F \longrightarrow \fin,\qquad \P(E) = (E\setminus\{0\})/\F^\times
\]
which associates to a \mod~$E$ over $\F$ the finite quotient space (set) $\P(E)$.

$b)$~It is important not to confuse the \mod~$\F^{(n)}$ with the induced module $\F\otimes_\B \B^{(n,1)}$ that can be realized
as the range $P(\F^n)$ of the projection $P$ as in  \eqref{projP} promoted to an element of $M_n(\F)$. There exists a unique map $\phi_n:\F\times\B^{(n,1)}\to\F^{(n)}=(\Z,\theta)^\vee$ which vanishes whenever one of the two arguments does so and   is defined as follows
$$
\phi_n:\F\times\B^{(n,1)}\to\F^{(n)}=(\Z,\theta)^\vee,\quad \phi_n(u^k,j)=\theta^k(j)=j+kn \qqq j\in \{0, \ldots,n-1\}.
$$
One has
 $\phi_n(x,i\vee j)=\phi_n(x,i)\vee \phi_n(x,j) $ $\forall i,j\in \{0, \ldots,n-1\}$ and $x\in \F$. Also $\phi_n(x\vee y,i)=\phi_n(x,i)\vee \phi_n(y,i) $ $\forall i\in \{0, \ldots,n-1\}$ and $x,y\in \F$. Moreover $\phi_n$ is $\F$-linear inasmuch as $\phi_n(\lambda x,j)=\lambda \phi_n(x,j)$ $\forall\lambda \in \F$.
Notice that there are more relations in $\F^{(n)}$ than those  holding in $\F\otimes_\B \B^{(n,1)}$: for example the relation~
$
\phi(1,i)\vee \phi(u,j)=\phi(u,j)\quad \forall i,j \in \{0, \ldots,n-1\}.
$
The \mod~$\F^{(n)}$, when viewed as a \mod~over $\B$, has free rank equal to $1$ ($\ruk_\B( \F^{(n)})=1$) while a similar conclusion fails to hold, as soon as $n>1$, for the \mod~$\B^{(n,1)}\otimes_\B \F$. 

As a set, $\F^{(n)}$ is the smash product $\F\wedge \B^{(n,1)}$ and its  additive structure  is given by the lexicographic order on the non-zero elements.  One needs to clarify in which sense this lexicographic smash product plays the role of the    tensor product $\otimes$ for $\B$-\mods~of free rank one.
}\end{rem}
The following table summarizes the geometric interpretation of the three categories $\Delta\subset \Lambda\subset \tilde \Lambda$ in terms of the geometric categories $\cP_\B\subset\cP^1_\F\subset\cP_\F$ ($\F=\Z_{max}$):

\bigskip

\begin{center}
\begin{tabular}{|c|c|}
\hline &   \\
Projective geometry $\cP_\B$ over $\B$ & Simplicial category $\Delta\sim\cP_\B$   \\
&  \\
\hline &   \\
Projective geometry $\cP_\F^1$ over $\F$ (linear)&  Cyclic category $\Lambda\sim \cP_\F^1$ 
 \\
&  \\
\hline &   \\
Projective geometry $\cP_\F$ over $\F$ (semilinear) &  Epicyclic category $\tilde\Lambda \sim \cP_\F$   \\
&  \\
\hline 
\end{tabular}
\end{center}


\bigskip

\subsection{The perfection of $\Z_{max}$ }

Let $\bff=\Q_{\rm max}:=(\Q\cup\{-\infty\},\text{max},+)$ be the sub-semifield of $\rmax$ (written multiplicatively) containing the (sub-)semifield $\F=\Z_{\max}$ (generated by the element $u>1$ of  $\rmax$)
as well as all rational powers $u^\alpha$, $\alpha \in \Q$. 

Recall that a semifield $K$ of characteristic one is called perfect when the 
map $x\mapsto x^n$ is surjective $\forall n\in \N^\times$. The map $\ff_n: x\mapsto x^n$ defines an automorphism of $K$ and one obtains an action $\ff$ of the multiplicative group $\Q_+^\star$ on $K$ such that $\ff_\alpha=\ff_n\circ \ff_m^{-1}$ for $\alpha=n/m$.

The following statement summarizes the main properties of $\bff=\Q_{max}$

 \begin{prop}
$(i)$~The semifield $\bff$ is perfect, contains $\F$ and for any perfect semifield $K\supset \F$
 one derives a canonical homomorphism $\bff\to K$ extending  the
inclusion $\F\subset K$.

$(ii)$~Any finitely generated sub-semifield of $\bff$ containing $\F$ is the inverse image $ \ff_n^{-1}(\F)=\ff_{1/n}(\F)$ of
$\F$ for some integer $n\in\N$. For $m,n\in\N$ one has $\ff_n^{-1}(\F)\subset \ff_m^{-1}(\F)$ if and only if $n\vert m$.

$(iii)$~The intersection of the semifields $ \ff_m(\F)$ $\forall m\in\N$ is the semifield $\B$.

\end{prop}
\proof $(i)$~A semifield $K$ of characteristic one is  perfect if and only if its multiplicative group $K^\times$ is uniquely divisible. This fact implies $(i)$ since $\Q\supset \Z$ is the uniquely divisible closure of $\Z$.

$(ii)$~Let $K\subset \bff$ be a  finitely generated sub-semifield of $\bff$ containing $\F$. Since sums $s=\sum b_j$ of elements of $\bff$ give one of the $b_j$, the multiplicative subgroup $K^\times\subset K$ is a finitely generated subgroup of $({\bff})^\times\sim \Q$ thus is of the form $\frac{1}{n}\Z\subset \Q$. This implies $(ii)$.

$(iii)$~Immediate. \endproof

\begin{rem}{\rm
In a semifield $K$ of characteristic one,  $x=1$ is  the only  solution of the equation $x^n=1$ since the endomorphism $\ff_n$ is injective. It follows that for any proper extension $K\supsetneq \bff$ the group 
$K^\times/{\bff}^\times$ is infinite and torsion free. Indeed $K^\times$ is torsion free and the group $({\bff})^\times\sim \Q$ is divisible, thus $K^\times/{\bff}^\times$ is also torsion free. }
\end{rem}

\section{Duality}\label{dual}

One key property of the cyclic category $\Lambda$ is that it is anti-isomorphic to itself, \ie one has a contravariant functor 
$\Lambda\longrightarrow\Lambda$, $f\mapsto f^t$, that determines an isomorphism $\Lambda\sim\Lambda^{\rm op}$. In this section we show that this duality corresponds to transposition in the framework of projective geometry in characteristic one as
developed in \S\ref{sectF1}.

\subsection{Self-duality of $\Lambda$}

\begin{prop}\label{propdual11/2}
$(i)$~Let $a,b$ be $>0$ integers and $f:\Z\to \Z$ be a non-decreasing map such that $f(x+a)=f(x)+b$ for all $x\in \Z$. Then there exists a unique map $f^t:\Z\to \Z$ such that 
\begin{equation}\label{dual0}
 f(x)\geq y\iff x\geq f^t(y) \qqq x,y\in \Z.
\end{equation}
Moreover $f^t$ is non-decreasing and fulfills $f^t(x+b)=f^t(x)+a$ for all $x\in \Z$.

$(ii)$~Let $h\in \Hom_\Lambda([n],[m])$,  then the class of $f^t$  is independent of the choice of $f\in h$ and defines an element $h^t\in \Hom_\Lambda([m],[n])$.

$(iii)$~The association $h\mapsto h^t$ defines a contravariant functor  $\Lambda\longrightarrow\Lambda$.
\end{prop}
\proof $(i)$~The inverse image $f^{-1}(I)$ of an interval $I\subset \Z$ of length $b$ is an interval $J$ of length $a$, since the translates $I+kb$ form a partition of $\Z$ and $f^{-1}(I+kb)=J+ka$. It follows that $f^{-1}([y,\infty))$ is of the form $[f^t(y),\infty)$
for a uniquely defined map $f^t:\Z\to \Z$. For $y\leq y'$ one has $f^{-1}([y',\infty))\subset f^{-1}([y,\infty))$ and thus $f^t(y')\geq f^t(y)$ so that $f$ is non-decreasing. Moreover the equality $f^{-1}(I+kb)=J+ka$ shows that $f^t(x+b)=f^t(x)+a$ $\forall x\in \Z$.\newline
$(ii)$~Let $f: \Z\to \Z$ be as in $(i)$. Let $k\in \Z$ and $f'$ be given by $f'(x)=f(x)-kb$ $\forall x\in \Z$. Then one has for any interval $I\subset \Z$ of length $b$, $f'^{-1}(I)=f^{-1}(I+kb)=J+ka$. This shows that $f'^t(y)=f^t(y)+ka$ $\forall y\in \Z$. Taking $a=n+1$, $b=m+1$ and using $(i)$, one obtains the required statement.\newline
$(iii)$~follows from the equalities
$
(f\circ g)^{-1}[y,\infty)=g^{-1}(f^{-1}[y,\infty))=g^{-1}[f^t(y),\infty)=[g^t(f^t(y)),\infty)
$
which show that $(f\circ g)^t=g^t\circ f^t$. \endproof

\begin{prop}
$(i)$~The covariant functor  $\Lambda\longrightarrow \Lambda$ which is the square of $h\mapsto h^t$ is equivalent to the identity by the natural transformation implemented by the map $[n]\mapsto \tau\in \Aut_\Lambda([n])$, $\tau(x)=x-1$ $\forall x\in \Z$.

$(ii)$~Let ${\rm Ad}(\tau)\in \Aut(\Lambda)$ be the inner automorphism defined by  the map $[n]\mapsto \tau\in \Aut_\Lambda([n])$ as in $(i)$. The action of the group $\Z$ on $\Lambda$: $j\mapsto {\rm Ad}(\tau)^j$  extends to a continuous action of the profinite completion $\hat\Z$: $\alpha:\hat \Z\to \Aut(\Lambda)$.

$(iii)$~An automorphism $\gamma\in  \Aut(\Lambda)$ belongs to $\alpha(\hat \Z)$ if and only if it is inner and its extension to 
$\Arc\triangleleft \N$ fixes the $\psi_k$ of \eqref{kmorpa}.
\end{prop}

\proof $(i)$~Let $h\in \Hom_\Lambda([n],[m])$ and $f\in h$, then $f^t\in h^t$ so that \eqref{dual0} holds. Using 
\eqref{dual0} one obtains 
$
f(x)< y\iff x< f^t(y)~  \forall x,y\in \Z
$
or equivalently
$
 f^t(y)\geq x+1\iff y\geq f(x)+1~ \forall x,y\in \Z
$.
One also has 
$
 f^t(y)\geq x+1\iff y\geq (f^t)^t(x+1)~ \forall x,y\in \Z
$
which thus gives the equality $ (f^t)^t(x+1)=f(x)+1$ $\forall x\in \Z$.\newline
$(ii)$~Let $z=(z_a)_{a\in \N}\in \prod_\N \Z/a\Z$, then,  the map $[n]\mapsto \tau^{z_{n+1}}\in \Aut_\Lambda([n])$ implements an  inner automorphism $\beta(z)\in \Aut(\Lambda)$ and the map $\beta: \prod_\N \Z/a\Z\to \Aut(\Lambda)$ is a continuous group homomorphism. Composing $\beta$ with the natural inclusion $\hat \Z\subset \prod_\N \Z/a\Z$, one obtains the required continuous action $\alpha:\hat \Z\to \Aut(\Lambda)$.\newline
$(iii)$~Let  $\gamma\in  \Aut(\Lambda)$ be an inner automorphism. Then we claim that there exists a unique $z=(z_a)_{a\in \N}\in \prod_\N \Z/a\Z$ such that $\gamma=\beta(z)$. Indeed, since every element of $\Aut_\Lambda([n])$ is a power $\tau^{z_{n+1}}$  one gets the existence of $z$, the uniqueness follows since the action of $\gamma$ on $\Hom_\Lambda([0],[n])$ uniquely determines $z_{n+1}$ modulo $n+1$. For a pair $n,k\in \N$ one has  $\psi_k\in  \Hom_{\Arc\triangleleft \N}(\Psi_k(\hat n),\hat n)$ and 
$
\tau^a \psi_k= \psi_k \tau^b \iff b\equiv a \ {\rm modulo} \ n.
$
Thus the extension of the inner automorphism $\gamma=\beta(z)$ to $\Arc\triangleleft \N$    fixes the $\psi_k$ if and only if 
\begin{equation}\label{projlim}
z_b\equiv z_a \ {\rm modulo} \ a  \qqq b=ka
\end{equation}
In turns \eqref{projlim} characterizes the elements of the projective limit $\hat \Z=\varprojlim \Z/a\Z$, thus $(iii)$ holds. \endproof

\subsection{Duality and transposition for \mods}

Next we describe the relation between the contravariant functor $\Lambda\longrightarrow\Lambda$, $f\mapsto f^t$   and the transposition of morphisms in linear algebra. Transposition is determined  in a  precise form by implementing the duality for $\B$-\mods~$E$ with $\ruk E=1$. Recall that for any such $\B$-\mod~$E$ the relation $x\leq y\iff x\vee y=y$ is a total order on $E$.

\begin{prop}\label{propdual}
$(i)$~Let $E$ be a $\B$-\mod~with $\ruk E=1$. Then for any $y\in E$ the following formula defines a linear form $\ell_y\in \Hom_\B(E,\B)$
\begin{equation}\label{idem3}
\ell_y(x)=< x , y>_\B:= \left\{
             \begin{array}{ll}
               0, & \hbox{if}\ x\leq y \\
               1, & \hbox{if}\ x> y.
             \end{array}\right.
\end{equation}
$(ii)$~For $z,t\in E$ set $ z\wedge t:= \inf(z,t)$. The pairing \eqref{idem3} satisfies the following bilinearity property
\begin{equation}\label{idem4}
 < x \vee y , z \wedge t>_\B= <x,z>_\B+<y,z>_\B+<x,t>_\B+<y,t>_\B
\end{equation}
where $+$ denotes the idempotent addition in $\B$.

$(iii)$~Let $E$ be a $\B$-\mod~with $\ruk E=1$ and $\rk E<\infty$. Let $E^*$ be the set $E$ endowed with the binary operation $\wedge$ as in $(ii)$. Then $E^*$ is a $\B$-\mod~with $\ruk E^*=1$, $\rk E^*=\rk E$. Moreover, the map 
$y\mapsto \ell_y$ defines a $\B$-linear isomorphism $E^*\sim  \Hom_\B(E,\B)$.
\end{prop}
\proof $(i)$~One has $\ell_y(0)=0$ since $0\leq y$, $\forall y\in E$. Moreover for any two elements $x,x'\in E$ the following equality holds $ < x \vee x' , y>_\B= <x,y>_\B+<x',y>_\B$ since one of these elements is $>y$ if and only if the largest of the two is $>y$.\newline
$(ii)$~For $y,z\in E$, with  $y\wedge z:= \inf(y,z)$ one has: 
$
< x, y\wedge z>_\B=<x,y>_\B+<x,z>_\B, \  \forall  x\in E.
$
Indeed: $\inf(y,z)<x$ if and only if $y<x$ or $z<x$.\newline
$(iii)$~We can view $E$ as a finite totally ordered set, then $E^*$ is the same set but endowed with the opposite ordering
so that the largest element of $E$ is the smallest in $E^*$ \ie the $0$-element for $E^*$. It follows that $\ruk E^*=1$, $\rk E^*=\rk E$. The map $E^*\to  \Hom_\B(E,\B)$ 
$y\mapsto \ell_y$ is $\B$-linear by \eqref{idem4}. It is injective since  $ \ell_y^{-1}(\{0\})=[0,y]$. We show that it is also surjective. Let $L\in \Hom_\B(E,\B)$, then $L(0)=0$ and 
$L$ is a non-decreasing so $L^{-1}(\{0\})=[0,y]$ for some $y\in E$, thus $L=\ell_y$.\endproof

Proposition \ref{propdual} shows that the duality of $\B$-\mods~$E$ with $\ruk E=1$ and $\rk E <\infty$ behaves similarly to the duality holding for finite 
dimensional vector spaces over fields and it produces in particular the transposition of linear maps defined as follows. Let $E^*=\Hom_\B(E,\B)$, $F^*=\Hom_\B(F,\B)$
\begin{equation*}
\Hom_\B(E,F) \ni f\mapsto f^*\in \Hom_\B(F^*,E^*), \quad    f^*(L)=L\circ f\qqq L\in \Hom_\B(F,\B).
\end{equation*}

As a corollary of Proposition \ref{propfin}, one derives that the simplicial category $\Delta$ is canonically isomorphic to the full subcategory of the category of $\B$-\mods~whose objects are the $\B^{(n,1)}$ for $n\geq 1$ and the morphisms are the linear maps $f$ such that $f^{-1}(\{0\})=\{0\}$. 

Although one has a canonical isomorphism $(\B^{(n,1)})^*\sim \B^{(n,1)}$, the condition $f^{-1}(\{0\})=\{0\}$ is  {\em not preserved} by transposition. In the following discussion we shall show that the transposed of the above condition is understood within the category $\cJ$ of {\em intervals}  (\ie totally ordered sets) $I$ with a smallest element $b_I\in I$ ($b_I\leq a$, $\forall a\in I$) and a largest element $t_I\in I$: \cf~\cite{MM} VIII.8.

The morphisms of the category $\cJ$ are
\begin{equation*}
   \Hom_\geq(I,J)=\{f:I\to J\mid x\leq y\implies f(x)\leq f(y), \ f(b_I)=b_J, \ f(t_I)=t_J\}
\end{equation*}
namely non-decreasing maps preserving the two end points. Notice that given a $\B$-\mod~$E$ with $\ruk(E) =1$ and $\rk(E) <\infty$, the underlying ordered set $E_\leq$ is an interval.
\begin{prop}\label{propdual1}
Let $f\in \Hom_\B(E,F)$, with $E,F$ $\B$-\mods~of finite rank and free rank  $1$. Then \newline
$(i)$~$f^{-1}(\{0\})=\{0\}$ if and only if $f^*\in \Hom_\leq(F_\leq^*,E_\leq^*)$.\newline
$(ii)$~$f\in \Hom_\leq(E_\leq,F_\leq)$ if and only if  $f^*\in \Hom_\B(F^*,E^*)$ fulfills $(f^*)^{-1}(\{0\})=\{0\}$.\newline
$(iii)$~The transposition of maps $f \mapsto f^*$ determines an isomorphism of  $\Delta^{\rm op}$ with the full subcategory of $\cJ$ defined by the intervals of the form $
n^*:=\{0,1,\ldots,n+1\}$, for $n\geq 0$. 
\end{prop}
\proof $(i)$~Let $E$ be a $\B$-\mod~with $\ruk E=1$ and $\rk E<\infty$. Then $\Hom_\B(E,\B)$ is an interval whose largest element is the linear form $\tau_E$: $\tau_E(x)=1$ $\forall x\in E$, $x\neq 0$. For $f\in \Hom_\B(E,F)$ as in $(i)$, one has $f^{-1}(\{0\})=\{0\}$ if and only if $\tau_F\circ f=\tau_E$. The smallest element of $\Hom_\B(F,\B)$ is the linear form $0$ which is automatically preserved by composition with any $f\in \Hom_\B(E,F)$.\newline
$(ii)$~One has $f\in \Hom_\leq(E_\leq,F_\leq)$ if and only if $f(t_E)=t_F$ where $t_E$ (resp. $t_F$) is the largest element of $E$. This holds if and only if  $t_E \notin f^{-1}([0,y])$ $\forall y<t_F$ \ie if and only if $(f^*)^{-1}(\{0\})=\{0\}$.\newline
$(iii)$~For $n\geq 0$, the dual \mod~$(\B^{(n+1,1)})^*$  is an interval of cardinality $n+2$, hence coincides with $n^*$. Transposition determines a contravariant functor $\cP_\B\longrightarrow\cJ$.
\endproof
In terms of the isomorphism $E^*\to  \Hom_\B(E,\B)$ 
$y\mapsto \ell_y$ of Proposition \ref{propdual}, the transposed $f^*$ of $f\in \Hom_\B(E,F)$ replaces $\ell_y$ by $\ell_y\circ f$, for $y\in F$, hence is defined by the equation
 \begin{equation}\label{dual1}
<f(x), y>_\B= <x,f^*(y)>_\B \qqq x\in E, y \in F^*.
\end{equation}

Using the above notations we obtain the following description  for the basic equation \eqref{dual0}
 \begin{equation}\label{dual2}
<y, f(x)>_\B= <f^t(y),x>_\B \qqq x,y\in \Z.
\end{equation}
This shows  that once interpreted in the framework of characteristic one, the contravariant functor $\Lambda\longrightarrow\Lambda$, $f\mapsto f^t$, is simply inverse transposition.

\begin{rem}
{\rm At first, it might seem puzzling that the transposition $f\mapsto f^*$ fulfills $(f^*)^*=f$ while the map $f\mapsto f^t$ 
of Proposition \ref{propdual1} is not involutive. The reason for this behavior is that for a finite totally ordered set $E$ viewed as a $\B$-\mod~the dual $E^*$ is the same set  but endowed with the {\em opposite} order: $x\leq^*y \iff y\leq x$. By applying  \eqref{dual1} one derives
$
f^*(z)\leq^*t \iff z \leq^*(f^*)^*(t),~ \forall z\in F^*, t\in E^*
$
which is equivalent to 
$
t\leq f^*(z) \iff (f^*)^*(t)\leq z.
$
This shows, using \eqref{dual1} that $(f^*)^*=f$. When considering the  map $f\mapsto f^t$ 
of Proposition \ref{propdual1}, one applies the same formula twice, while instead 
 taking into account the {\em opposite} order would provide the inverse of the map $f\mapsto f^t$. Since the negation of 
$x\leq y$ is $x>y$ \ie $y+1\leq x$ for the ordered set $\Z$, the translation of $1$ pops-up and conjugates the map $f\mapsto f^t$ with its inverse.
}\end{rem}

Next we develop the above duality directly at the level of \mods~over $\F=\Z_{max}$. Let $(X,\theta)$ be an archimedean set and let $E=(X,\theta)^\vee$ be the associated \mod~over $\F$ as in \eqref{vect1} and \eqref{vect2}. The archimedean property ensures that the following  pairing is well defined with values in $\F$:
 \begin{equation}\label{pairingF}
<x, y>_\F:=\inf \{v\in \F\mid x \leq vy\}\qqq x,y \in E, y\neq 0.
\end{equation}

\begin{prop}\label{propdual2}

$(i)$~Let $(X,\theta)$ be an archimedean set and $E=(X,\theta)^\vee$ the associated \mod~over $\F$ (\cf~Proposition~\ref{my check}). Let $E^*=(X',\theta^{-1})^\vee$ where $X'$ is the set $X$ endowed with the opposite order. Then \eqref{pairingF} defines a bilinear 
pairing  $E\times E^*\to\F$.\newline
$(ii)$~Let $X=\Z$, $\theta(x)=x+n$, and $E$, $E^*$ the associated \mods~over $\F$ as in $(i)$. Then \eqref{pairingF} determines an isomorphism 
 \begin{equation*}
\ell: E^*\stackrel{\sim}{\to} \Hom_\F(E,\F), \  \  \ell_y(x):=<x, y>_\F \qqq x \in E, y\in E^*.
\end{equation*}
$(iii)$~The contravariant functor $\Lambda\longrightarrow\Lambda$ $f\mapsto f^t$  is the  inverse of transposition $f\mapsto f^*$:
 \begin{equation}\label{pairingF2}
\Hom_\F(E,F)\ni f \mapsto f^*\in \Hom_\F(F^*,E^*) \  \   <f(x),y>_\F= <x,f^*(y)>_\F \qqq x\in E, y\in F^*.
\end{equation}
\end{prop}
\proof $(i)$~For $y\in X$ the archimedean property ensures that the set $ \{v\in \F\mid x \leq vy\}$ is non empty $\forall x\in E$. Thus \eqref{pairingF} is well-defined and gives $0\in \F$ only for $x=0\in E$. For $0\neq \lambda\in \F$ one has 
$
\{v\in \F\mid \lambda x \leq vy\}= \lambda \{v\in \F\mid x \leq vy\}
$
which shows that $<\lambda x, y>_\F=\lambda<x, y>_\F$. For $x\leq x'$ one has 
$
\{v\in \F\mid x' \leq vy\}\subset \{v\in \F\mid x \leq vy\}
$,
thus $<x,y>_\F \leq <x',y>_\F$. This shows that the map $x\mapsto <x,y>_\F\in \F$ is $\F$-linear. For $x,y\in E$ and $y\neq 0$ one has for 
$0\neq \lambda \in \F$: $< x, \lambda y>_\F=\lambda^{-1}<x, y>_\F$ which corresponds to the $\F$-linearity in $E^*$. Also one has $<x,y'>_\F\leq <x,y>_\F$ for $y \leq y'$ for $y,y'\in E$ this corresponds to the linearity in $E^*$. Thus the pairing \eqref{pairingF} is bilinear. \newline
$(ii)$~It follows from $(i)$ that the map $\ell: E^*\to \Hom_\F(E,\F)$ is well defined and linear. Let $L\in \Hom_\F(E,\F)$, we show that there exists a unique $y\in E^*$ with $\ell_y=L$. This holds for $L=0$ thus we can assume that $L(x_0)\neq 0$ for some $x_0\in E$. The archimedean property implies that the kernel of $L$, \ie $\{x\in E\mid L(x)=0\}$ is reduced to $\{0\}$ since for any $x\in E$ there exists  $\lambda\in \F$ with $x\geq \lambda x_0$ so that $L(x)\geq \lambda L(x_0)\neq 0$. Replacing $L$ by a multiple $\mu L$ for some $\mu \neq 0$ we can thus assume that $L$ corresponds, at the level of archimedean sets, to a non-decreasing map $f:\Z \to \Z$ with $f(x+n)=f(x)+1$ $\forall x\in \Z$ and that $f(0)=0$. On then has $f(n)=1$ and $f$ is uniquely determined by the element $y\in [0,n-1]$ such that $f(x)=0\iff x\leq y$. We show that $L=\ell_y$, \ie that for any $x\in \Z$, $f(x)=g(x)$ where $g(x)$ is the smallest integer $k\in \Z$ such that $x\leq y+kn$. Since $g(x+n)=g(x)+1$ $\forall x\in \Z$, it is enough to prove that $f(x)=g(x)$ for $x\in [0,n-1]$. For $x\in [0,y]$ one gets $g(x)=0$ since $y-n<x\leq y$. For $x\in [0,n-1]$, $x>y$, one has $g(x)=1$, since $y<x\leq y+n$. This shows that $L=\ell_y$ for a unique $y\in \Z$ and thus that the map $\ell:E^*\to \Hom_\F(E,\F)$ is bijective.\newline
$(iii)$~The equality 
\begin{equation*}
\Re_\B(x)= \left\{
             \begin{array}{ll}
               0, & \hbox{if}\ x\leq 1 \\
               1, & \hbox{if}\ x> 1
             \end{array}\right. \qqq x\in \F
\end{equation*}
defines a linear form $\Re_\B:\F\to \B$. Moreover, for $(X,\theta)$  an archimedean set and $E=(X,\theta)^\vee$ the associated \mod~over $\F$, one has 
\begin{equation*}
<x,y>_\B=\Re_\B(<x,y>_\F)\qqq x,y\in E, y\neq 0.
\end{equation*}
Indeed $x\leq y$ if and only if $1\in \{v\in \F\mid x \leq vy\} $.
By construction, the transposition $\Hom_\F(E,F)\ni f \mapsto f^*\in \Hom_\F(F^*,E^*)$ fulfills \eqref{pairingF2}. Applying $\Re_\B$ to both sides of \eqref{pairingF2} one obtains 
 \begin{equation*}
  <f(x),y>_\B= <x,f^*(y)>_\B \qqq x\in E, y\in F^*
\end{equation*}
and this shows, using \eqref{dual2}, that $(f^*)^t=f$.\endproof

\subsection{Lift of permutations and cyclic descent number}

Following a traditional point of view  the symmetric group $S_n$ is interpreted as the limit for $q\to 1$ of the general linear group ${\rm GL}_n(\F_q)$ over a finite field $\F_q$. Notice that in the limit  the cardinality of the projective space $\P^{n-1}(\F_q)$ becomes $n$.  In Proposition \ref{propsemi2} we have proven
that $\cF_1$ extends to  a functor $\cF:\Arc\triangleleft \N\longrightarrow \Se$. When interpreted in terms of geometry over the semifield $\F=\Z_{max}$, this functor associates to a \mod~$E$ over $\F$ the quotient set $(E\setminus \{0\})/\F^\times$. Now we restrict  this functor to the epicyclic category $\tilde \Lambda$ \ie to the \mods~$\F^{(n)}$ obtained from the one dimensional free vector space $\F$ by restriction of scalars as explained in Theorem \ref{thmepi}. Note that since $\End(\F)\sim \N^\times$, the functor ${\rm Mod}$ restricts to a functor ${\rm Mod}:\cP_\F\longrightarrow \N^\times$ where $\N^\times$ is viewed as a category with a single object.

\begin{defn} (\cf\cite{Loday})~
Let $\sigma\in \Hom_\fin([n],[m])$ be a map of sets, then the cyclic descent number 
of $\sigma$ is defined to be
\begin{equation*}
{\rm cdesc}(\sigma)=\#\{j\in \{0,1,\ldots ,n\}\mid \sigma(j+1)<\sigma(j)\}
\end{equation*}
\end{defn}
where we identify $n+1\sim 0$.

The following result provides a geometric interpretation of the cyclic descent number of an arbitrary permutation as the measure of its semilinearity
\begin{prop}\label{propdescent}
$(i)$~The functor $\P:\cP_\F\longrightarrow \fin$  is full.\newline
$(ii)$~Let $\sigma\in \Hom_\fin([n],[m])$. Then ${\rm cdesc}(\sigma)$  is the smallest integer $k$ such that 
there exists $
 f\in \Hom_{\cP_\F}(\F^{(n)},\F^{(m)})$,  ${\rm Mod}(f)=k$ with $ \P(f)=\sigma$.
\newline
$(iii)$~Let $\sigma\in \Hom_\fin([n],[m])$ with ${\rm cdesc}(\sigma)=k$,  then 
there exists a unique $
 f\in \Hom_{\cP_\F}(\F^{(n)},\F^{(m)})$,  ${\rm Mod}(f)=k$ such that $ \P(f)=\sigma$.
\end{prop}
\proof It is enough to prove $(ii)$ and $(iii)$.   Let $U_k$ be the set of  $
 f\in \Hom_{\Arc\triangleleft \N}([n],[m])$,  $\rho(f)=k$ (\cf~Proposition~\ref{propsemi} (ii)) such that $ \cF (f)=\sigma$. Using Theorem \ref{thmepi} it is enough to show that $U_k= \emptyset$ for $k<{\rm cdesc}(\sigma)$ and that if $k={\rm cdesc}(\sigma)$ then $U_k$ contains a single element. By construction (\cf~Definition \ref{last2}) the elements of $U_k$ are equivalence classes of non-decreasing maps, modulo the addition of a constant multiple of $m+1$, $f:\Z\to \Z$ such that 

$a)$~$f(x+(n+1))=f(x)+ k(m+1)$ $\forall x \in \Z$

$b)$~$f(x)\in \sigma(x)+(m+1)\Z$~ $\forall x\in\{0,1,...,n,n+1\}$, $\sigma(n+1):=\sigma(0)$.

In each equivalence class there is a unique representative $f$ such that $f(0)=\sigma(0)$: in the following we assume this normalization condition. We let
\begin{equation}\label{sem1}
c(x)=\#\{j\in \{0,1,\ldots ,x-1\}\mid \sigma(j+1)<\sigma(j)\}, \ \ \forall x\in\{0,1,\ldots ,n,n+1\}
\end{equation} 
where by convention we set $\sigma(n+1):=\sigma(0)$ so that $c(n+1)={\rm cdesc}(\sigma)$.
For $x\in\{0,1,...,n,n+1\}$, let $b(x)\in \Z$ such that $f(x)=\sigma(x) +(m+1) b(x)$. One has $b(0)=0$ and since $f(1)\geq f(0)$ we get $b(1)\geq c(1)$. More generally,  for $j\in\{0,1,\ldots ,n\}$, one obtains
$$
f(j+1)\geq f(j)\implies b(j+1)-b(j)\geq c(j+1)-c(j).
$$
Indeed,  one has $(m+1)(b(j+1)-b(j))\geq (\sigma(j)-\sigma(j+1))$ and this implies
 $$\Z \ni b(j+1)-b(j)\geq \frac{\sigma(j)-\sigma(j+1)}{m+1}>-1.$$
If $\sigma(j+1)<\sigma(j)$ then $b(j+1)-b(j)\geq 1=c(j+1)-c(j)$.

The inequalities $b(j+1)-b(j)\geq c(j+1)-c(j)$ for $j\in\{0,1,\ldots ,n\}$ together with $b(0)=c(0)=0$ show that $b(n+1)\geq c(n+1)$. By $a)$ one has $b(n+1)=k$ while $c(n+1)={\rm cdesc}(\sigma)$. Thus $U_k= \emptyset$ for $k<{\rm cdesc}(\sigma)$. Moreover for $k={\rm cdesc}(\sigma)$, all the inequalities $b(j+1)-b(j)\geq c(j+1)-c(j)$ for $j\in\{0,1,\ldots ,n\}$ become equalities and 
 we obtain
\begin{equation}\label{sem1bis}
f(x)=\sigma(x)+ (m+1) c(x), \  \forall x\in\{0,1,\ldots ,n\}
\end{equation}
which provides the required uniqueness, thus $U_k$ contains at most one element. Moreover one easily checks that the function $f$ defined by  \eqref{sem1bis} on $\{0,1,\ldots ,n\}$ and extended by periodicity using $a)$ is non-decreasing and belongs to $U_k$. \endproof

\section{Extension of scalars $-\tensb$ to hyperfields}\label{hypersec}
The algebraic  constructions discussed in the earlier sections for semifields are in fact the ``positive part''
of a general picture that one can elaborate in terms of hyperfields. The development of this translation into the framework of hyperstructures allows one to link with the results of \cite{CC5} where it was  shown that by implementing the theory of hyperstructures  one can parallel successfully 
Fontaine's $p$-adic arithmetic theory of ``perfection'' and subsequent Witt extension by combining a process of dequantization (to characteristic one) and a consecutive Witt construction (to characteristic zero).
 It turns out that the \mods~implemented over the semifields $\B$ and $\F$ of last sections
fulfill precisely the property \eqref{share}  below that allows one to apply the symmetrization process introduced in \cite{Henry}. This procedure associates to a commutative mono\"{\i}d $M$ such that
\begin{equation}\label{share}
\forall x,y,u,v \in M, \  \ x+y = u+v\implies \exists z \in M, \  \
\begin{cases} x+z=u, & z+v=y, 
\\
\text{or}\  \ x=u+z, & v=z+y.
\end{cases}
\end{equation}
a hypergroup $s(M)$ which is the {\em universal solution} to the embedding of $M$ into a hypergroup.

It is shown in \cite{Henry} that the condition \eqref{share} is equivalent to the existence of a common refinement of any two decompositions of an element of $M$ as a sum. 

Let  $E$ be a \mod~over $\B$, then using the results of \opcit one shows that $E$, as a mono\"{\i}d, fulfills \eqref{share} if and only if $\ruk E=1$.  Moreover, in  \cite{Henry} is also proven that  the hypergroup $s(E)$ which is the   universal solution to the embedding of $E$ into a hypergroup  coincides with the tensor product $E\tensb$ which we now describe in details. 
Let  $E$ be a \mod~over $\B$ such that $\ruk E=1$. We denote by $E\tensb$ the quotient of $E\times \{\pm 1\}$ by the equivalence relation that identifies $(0,-1)\sim (0,1)$. We use the notation $\pm x$ for the elements of $E\times \{\pm 1\}$ and denote by $x\mapsto \vert x\vert$ the projection from $E\tensb$ to $E$, so that $\vert \pm x\vert :=x$. We endow   $E\tensb$ with the multivalued binary operation (here we use the total order $<$ of $E$)
\begin{equation}\label{taurgood}
    x \smile y=\left\{
                 \begin{array}{ll}
                  x, & \hbox{if $|x|>|y|$ or $x=y$;} \\
                   y, & \hbox{if $|x|<|y|$ or $x=y$;} \\
                    $[$-x,x$]$, & \hbox{if $y=-x$}
                 \end{array}
               \right.
\end{equation}
where the interval $[-x,x]$ is defined to be the set $\{z\in E\tensb\mid \vert z\vert \leq \vert x\vert \}$.

The following lemma lists several properties inherent to the above defined tensor product
\begin{lem}\label{hyper}
Let $E,F$ be \mods~over $\B$ such that $\ruk E=\ruk F=1$.

$(i)$~$E\tensb$ is a canonical hypergroup and a module over $\sign$ in the following sense
\begin{equation}\label{weak}
\lambda (v\smile v')=\lambda v \smile \lambda v', \ \ \ \ \ (\lambda+\lambda') v\subset \lambda v \smile \lambda' v\qqq \lambda, \lambda '\in \sign, \quad \forall v, v'\in E\tensb.
\end{equation}
$(ii)$~Let $f:E\to F$ be a morphism of $\B$-\mods. Then the following formula defines a morphism of $\sign$-modules:
\begin{equation*}
f\tensid : E\tensb \to F\tensb, \ \  f\tensid(\pm x):= \pm f(x) \qqq x\in E.
\end{equation*}
$(iii)$~Let $\epsilon: (E\setminus \{0\})\to \{\pm 1\}$ be a map of sets. Then the following defines an automorphism $\tilde \epsilon: E\tensb \to E\tensb$ (as a module over $\sign$)
\begin{equation*}
\tilde \epsilon(\pm x):= \pm \epsilon(x) x \qqq x\in E.
\end{equation*}
$(iv)$~Let $g\in \Hom_\sign(E\tensb, F\tensb)$ be a morphism of $\sign$-modules such that $g^{-1}(\{0\})=\{0\}$. Then there exists a unique pair $(f,\epsilon)$ of a morphism of $\B$-\mods~
$f:E\to F$ and a map of sets $\epsilon: (E\setminus \{0\})\to \{\pm 1\}$ such that $g=(f\tensid )\circ \tilde \epsilon$.
\end{lem}
\proof 
$(i)$~We recall (\cf\cite{CC4}) that the definition of a  canonical hypergroup $H$ requires that $H$ has a neutral element $0\in H$ (\ie an additive identity) and that the following axioms apply\vspace{.05in}

 $(1)$~$x+y=y+x,\qquad\forall x,y\in H$\vspace{.05in}

$(2)$~$(x+y)+z=x+(y+z),\qquad\forall x,y,z\in H$ \vspace{.05in}

$(3)$~$0+x=x= x+0,\qquad \forall x\in H$\vspace{.05in}

$(4)$~$\forall x  \in H~  \ \exists!~y(=-x)\in H\quad {\rm s.t.}\quad 0\in x+y$\vspace{.05in}

$(5)$~$x\in y+z~\Longrightarrow~ z\in x-y.$\vspace{.05in}

$(1)$, $(3)$ and $(4)$ are easy to verify for $E\tensb$. To check $(2)$ note that if among $\vert x\vert , \vert y\vert , \vert z\vert $ only one, say $\vert x\vert$, is strictly larger than the others, then both sides of $(2)$ give $x$. For $\vert y\vert \leq \vert x\vert $ one has $y \smile [-x,x]=[-x,x]$ since $y \smile [-x,x]$ contains  $[-y,y]$ and any $z$ such that $\vert y\vert <\vert z\vert\leq \vert x\vert $. It follows that $(2)$ holds when, among $\vert x\vert , \vert y\vert , \vert z\vert $ two are equal and strictly larger than the remaining one. Finally when $\vert x\vert =\vert y\vert = \vert z\vert $ both sides give $[-x,x]$ except when $x=y=z$ in which case they both give $x$.  The condition $(5)$ follows from the first $4$ since both sides are equivalent to $0\in -x+y+z $. This shows that $E\tensb$ is a canonical hypergroup. Moreover one has $ -x \smile -y= -(x \smile y)$ $\forall x,y\in E\tensb$ which gives the first equality of \eqref{weak}. One finally checks that the second inclusion holds and is in general strict for $\lambda'=-\lambda$. \newline
$(ii)$~Let $g : E\tensb \to F\tensb$ be defined by  $g(\pm x):= \pm f(x)$, $ \forall x\in E$. By construction $g(\lambda x)=\lambda g(x)$ $\forall\lambda \in \sign$ and $\forall x\in E\tensb$. It remains to check that $g$ is a morphism of hypergroups, \ie that $g(x \smile y) \subset g(x)\smile g(y)$ $\forall x,y\in E\tensb$. Since $f$ is non-decreasing $\vert g(x)\vert \leq \vert g(y)\vert $ if $\vert x\vert \leq \vert y\vert $. If $\vert x\vert < \vert y\vert $ $x \smile y=y$, $g(x \smile y) =g(y)\in g(x)\smile g(y)$ since $v\in u\smile v$ when $\vert u\vert \leq \vert v\vert$. The only remaining case to consider is when $y=-x$. One has  $g(x \smile -x) =g([-x,x])\subset [-g(x),g(x)]$.\newline
$(iii)$~By construction $\tilde \epsilon(\lambda x)=\lambda \tilde \epsilon(x)$ $\forall\lambda \in \sign$ and $\forall x\in E\tensb$. It remains to show that $\tilde \epsilon(x \smile y) = \tilde \epsilon(x)\smile \tilde \epsilon(y)$ $\forall x,y\in E\tensb$. In fact one has $\vert  \tilde \epsilon(x)\vert =\vert x\vert$ $\forall x\in E\tensb$. Thus \eqref{taurgood} shows the required equality when $\vert x\vert \neq \vert y\vert$ or when $x=y$. The only remaining case is $y=-x$, and in that case the equality follows from $\tilde \epsilon([-x,x])=[-x,x]$.\newline
$(iv)$~The map $f:E\to F$ is uniquely determined by $f(x):=\vert g(x)\vert$ $\forall x\in E$. Since $g^{-1}(\{0\})=\{0\}$, this determines uniquely the map $\epsilon: E\setminus \{0\}\to \{\pm 1\}$ such that $g=(f\tensid )\circ \tilde \epsilon$. It only remains to show that $f$
is non-decreasing \ie that
$x \leq y$ implies $\vert g(x)\vert\leq \vert g(y)\vert$. Assume $x < y$ and $\vert g(x)\vert> \vert g(y)\vert$. By hypothesis one has $g(x \smile y) \subset g(x)\smile g(y)$ $\forall x,y\in E\tensb$ and this contradicts  \eqref{taurgood} which gives 
$x\smile y=y$ and $g(x)\smile g(y)=g(x)$.
\endproof
 Applying Lemma \ref{hyper} to the  semifields $\F\subset \bff$ ($\F = \Z_{max}$) one obtains the corresponding hyperfields $\F\tensb\subset \bff\tensb$. The hyperfield $\bff\tensb$ is perfect and coincides with the perfection of $\F\tensb$. Using this functorial construction, one can recast the results of the previous sections in terms of hyperfields. 
 
 Finally, notice that as a set,  the projective space $\P(E)$ remains unchanged after shifting to the framework of hyperfields since the multiplicative group  of \eg $\F\tensb$ is simply the product of $\F^\times$ by the group of signs $\{\pm 1\}=\sign^\times$. Thus for a \mod~$E$ over $\F$ the following equality of sets holds
$$
\left((E\tensb)\setminus \{0\}\right)/(\F\tensb)^\times=\P(E)=(E\setminus \{0\})/\F^\times.
$$

\end{document}